\input ./style/arxiv-ba.cfg
\documentclass[ba,linksfromyear,preprint]{imsart}
\makeatletter
   \@ifpackageloaded{natbib}{}{\usepackage{natbib}}
\makeatother
\usepackage{algorithmic,url,color,placeins}

\pubyear{2015}
\volume{10}
\issue{4}
\firstpage{937}
\lastpage{963}
\doi{10.1214/14-BA930}

\startlocaldefs
\newcommand{\enquote}[1]{``#1''}

\newcommand{\bbr}{\mathbb{R}}
\newcommand{\approptoinn}[2]{\mathrel{\vcenter{
\offinterlineskip\halign{\hfil$##$\cr
#1\propto\cr\noalign{\kern2pt}#1\sim\cr\noalign{\kern-2pt}}}}}
\newcommand{\appropto}{\mathpalette\approptoinn\relax}
\newcommand{\given}{\mathrel{|}} 

\newtheorem{definition}{Definition}
\newtheorem{algorithm}{Algorithm}
\endlocaldefs

\begin{document}

\begin{frontmatter}
\title{Sequential Bayesian Model Selection of Regular Vine Copulas}
\runtitle{Sequential Bayesian Model Selection of Regular Vine Copulas}

\begin{aug}
\author[a]{\fnms{Lutz} \snm{Gruber}\corref{}\ead[label=e1]{gruber@ma.tum.de}}
\and
\author[b]{\fnms{Claudia} \snm{Czado}}

\runauthor{L. Gruber and C. Czado}

\address[a]{Center for Mathematics, Technical University of Munich, Munich, Germany, \printead{e1}}
\address[b]{Center for Mathematics, Technical University of Munich, Munich, Germany}

\end{aug}

%
\begin{abstract}
Regular vine copulas can describe a wider array of dependency patterns
than the multivariate Gaussian copula or the multivariate Student's t
copula. This paper presents two contributions related to model
selection of regular vine copulas. First, our pair copula family
selection procedure extends existing Bayesian family selection methods
by allowing pair families to be chosen from an arbitrary set of
candidate families. Second, our method represents the first Bayesian
model selection approach to include the regular vine density
construction in its scope of inference. The merits of our approach are
established in a simulation study that benchmarks against methods
suggested in current literature. A real data example about forecasting
of portfolio asset returns for risk measurement and investment
allocation illustrates the viability and relevance of the proposed scheme.
\end{abstract}

%
\begin{keyword}
\kwd{dependence models}
\kwd{graphical models}
\kwd{reversible jump MCMC}
\kwd{\break multivariate statistics}
\kwd{multivariate time series}
\kwd{portfolio risk forecasting}
\end{keyword}

\end{frontmatter}


\section{Introduction}
The use of copulas in statistics allows for the dependence of random
variables to be modeled separately from the marginal distributions.
This property makes copulas a very convenient tool to be used by
statistical modelers (\cite{nelsen2006,mcneilfreyembrechts2005,cookekurowicka2006,kjoe2010}).
While many classes of bivariate
copulas, also called pair copulas, are known (\cite{joe2001}), there is
only a very limited number of multivariate copulas available with a
closed-form analytical expression. Additionally, these only cover
limited patterns of dependence. Regular vine copulas provide a solution
to this problem by offering a construction method to design
multivariate copula densities as products of only bivariate
(conditional) copula densities (\cite{joe1996,cooke2001}).

A $d$-dimensional regular vine copula is set up in two steps. A
sequence of $d-1$ linked trees $\mathcal{V} = (T_1, \ldots,
T_{d-1})$, called the regular vine, functions as the building
plan for the pair copula construction. Each of the $d-j$ edges of tree
$T_j$, $j=1{:}(d-1)$, corresponds to a bivariate copula density that is
conditional on $j-1$ variables. A copula family is selected for each of
these (conditional) bivariate building blocks from a set of bivariate
(parametric) candidate families $\mathbf{B}$. The mapping of the pair
copula families to the regular vine is denoted by $\mathcal{B}_{\mathcal
{V}}$, and the parameters, which depend on the choice of the pair
copula families $\mathcal{B}_{\mathcal{V}}$, are denoted by $\boldsymbol
{\theta}_{\mathcal{V}}$.

Model selection of regular vine copulas is a difficult task, given that
there exist $\frac{d!}{2} \times2^{\binom{d-2}{2}}$ different
$d$-dimensional regular vine tree structures alone (\cite
{morales2011}). To obtain the number of possible regular vine copulas
on $d$ dimensions, the number of possible regular vine tree structures
must be multiplied by the number of possible combinations of pair
copula families, $|\mathbf{B}|^\frac{d(d-1)}{2}$. In higher dimensions,
there are too many possible models to analyze all of them to select the
suitable few. To reduce the complexity of model selection, \cite
{dissmann2010} suggested a tree-by-tree approach, which selects the
trees $T_1, \ldots, T_{d-1}$ of the regular vine $\mathcal{V}$
sequentially. We present a combination of \cite{dissmann2010}'s
tree-by-tree strategy of complexity reduction with a proper Bayesian
model selection scheme.

Existing research on Bayesian model selection of vine copulas is
restricted to two relatively small subclasses of regular vine copulas,
D-vine copulas and C-vine copulas. \cite{smithminalmeidaczado2009}
developed a model indicator-based approach to select between the
independence copula and one alternative copula family at the pair
copula level for D-vine copulas. \cite{minczado2010} and \cite
{minczado2011} discuss a more flexible approach to estimate the pair
copula families of D-vine copulas that is based on reversible jump
Markov chain Monte Carlo (MCMC) (\cite{green1995}). Our proposed method
provides two substantial contributions to Bayesian selection of vine
copulas: firstly, our method applies to the general class of regular
vine copulas and is capable of selecting the pair copula families
$\mathcal{B}_{\mathcal{V}}$ from an arbitrary set of candidate families
$\mathbf{B}$, which is beyond the scope of \cite
{smithminalmeidaczado2009, minczado2010, minczado2011} and contains
their selection procedures as special cases; secondly, our method is
the first Bayesian approach to estimating the regular vine $\mathcal
{V}$ of a regular vine copula along with the pair copula families
$\mathcal{B}_{\mathcal{V}}$. The latter innovation eliminates an
unrealistic assumption all previous Bayesian selection procedures rest
on---that the regular vine $\mathcal{V}$ is already known.


The remainder is organized as follows. Section \ref{ch:rvine-copulas}
formally introduces regular vine copulas to the extent required. In
Section \ref{ch:rvine-est}, we present our new approach to Bayesian
model selection for regular vine copulas using reversible jump Markov
chain Monte Carlo. Section \ref{ch:simstudy} presents the results of a
simulation study to help establish the validity of our model selection
algorithm. Section \ref{ch:case} presents a real data example using our
model selection procedure to improve forecasting of risk metrics of
financial portfolios. We conclude with further remarks in Section \ref
{ch:conclusions}.


\section{Regular Vine Copulas} \label{ch:rvine-copulas}
A copula describes a statistical model's dependence behavior separately
from its marginal distributions (\cite{sklar1959}). The copula
associated with a $d$-variate cumulative distribution function (cdf)
$F_{1{:}d}$ with univariate marginal cdfs $F_1, \ldots, F_d$ is a
multivariate distribution function $C{:}[0,1]^d \to
[0,1]$ with \textit{Uniform(0,1)} margins that satisfies
\[
F_{1{:}d}\left(\textbf{x}\right) = C\left(F_1\left(x_1\right), \ldots,
F_d\left(x_d\right)\right), \textbf{x} \in\bbr^d \text{.}\vadjust{\eject}
\]
%
\subsection{Model Formulation}
\begin{figure}
\includegraphics{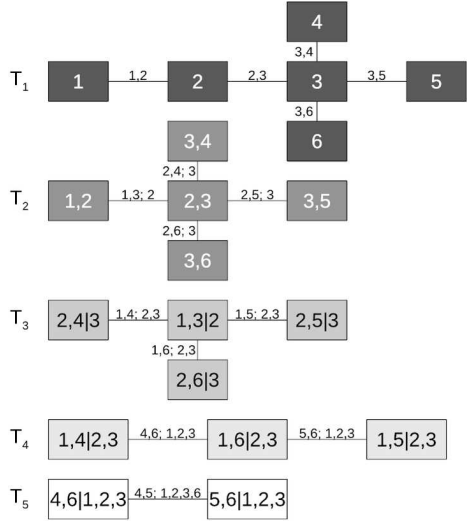}
\caption{Six-dimensional regular vine copula. This copula corresponds
to Scenario 1 of the simulation study of Section \ref{ch:simstudy}.}
\label{fig:ex-rv}
\end{figure}

\cite{joe1996} presented the first construction of a multivariate
copula using (conditional) bivariate copulas. \cite{cooke2001}
developed a more general construction method of multivariate densities
and introduced regular vines to organize different pair copula
constructions. The definitions and results stated in remainder of this
section follow \cite{cooke2001}, if not stated otherwise.
\begin{definition}[Regular vine tree sequence] \label{def:rvine}
A set of linked trees $\mathcal{V} = (T_1, T_2, \dots,
T_{d-1})$ is a regular vine on $d$ elements if
\begin{enumerate}
\item$T_1$ is a tree with nodes $N_1 = \{1, \dots, d\}$ and a set of
edges denoted by $E_1$.
\item For $k=2,\dots,d-1$, $T_k$ is a tree with nodes $N_k = E_{k-1}$
and edge set $E_k$.
\item For $k=2,\dots,d-1$, if $a=\{a_1, a_2\}$ and $b=\{b_1, b_2\}$ are
two nodes in $N_k$ connected by an edge, then exactly one of the $a_i$
equals one of the $b_i$ (Proximity condition).
\end{enumerate}
\end{definition}
Regular vines serve as the building plans for pair copula
constructions. When each edge of the regular vine is interpreted as a
(conditional) bivariate copula in the pair copula construction, the
resulting copula is a regular vine copula. We use the following
notation in the formal definition: the complete union $A_e$ of an edge
$e =\{a,b\} \in E_k$ in tree $T_k$ of a regular vine $\mathcal{V}$ is
defined by
\[
A_e = \left\{v \in N_1 \given\exists e_i \in E_i, i=1,\ldots,k-1, \text
{ such that } v \in e_1 \in\cdots\in e_{k-1} \in e\right\} \text{.}
\]
The conditioning set associated with edge $e =\{a,b\}$ is defined as
$D(e):=A_a \cap A_b$ and the conditioned sets associated with edge $e$
are defined as $i(e) :=A_a \setminus D(e)$ and $j(e) :=A_b \setminus
D(e)$, where $A \setminus B := A \cap B^c$ and $B^c$ is the complement
of $B$. The conditioned sets can be shown to be singletons (see, for
example, \cite{cookekurowicka2006}). In graphs we label an edge $e$ by
the derived quantities $i(e), j(e); D(e)$, which suggests a
probabilistic interpretation. Figure \ref{fig:ex-rv} shows a
6-dimensional regular vine to illustrate this notational convention. %
%


\begin{definition}[Regular vine copula] \label{def:rvine-cop}
Let $\mathcal{V} = (T_k = (N_k, E_k) \given k = 1, \ldots, d-1)$ be a
regular vine on $d$ elements. Let
\begin{equation*}
\mathcal{B}_k := (\mathcal{B}_e \given e \in E_k) \text{ and }
\boldsymbol{\theta}_k := (\boldsymbol{\theta}_e \given e \in E_k)
\end{equation*}
be the pair copula families and parameters of level $k$, where the
parameter vector $\boldsymbol{\theta}_e$ depends on the pair copula
family $\mathcal{B}_e$ of edge $e$. The pair copula families of all
levels $k=1{:}(d-1)$ are collected in $\mathcal{B}_{\mathcal{V}} :=
\mathcal{B}_{1{:}(d-1)} := (\mathcal{B}_1, \ldots, \mathcal{B}_{d-1})$.
The same notational convention extends to the parameters $\boldsymbol
{\theta}_{\mathcal{V}}$.

The regular vine copula $(\mathcal{V}, \mathcal{B}_{\mathcal{V}},
\boldsymbol{\theta}_{\mathcal{V}})$ has the density function
%
\begin{equation} \label{eq:rvine-pdf}
c(\mathbf{u}; \mathcal{V}, \mathcal{B}_{\mathcal{V}}, \boldsymbol{\theta
}_{\mathcal{V}}) = \prod_{T_k \in\mathcal{V}} \prod_{e \in E_k}
c_{\mathcal{B}_e}\left(u_{i(e) \given D(e)}, u_{j(e) \given D(e)};
\boldsymbol{\theta}_e\right) \text{,}
\end{equation}
where $c_{\mathcal{B}_e}(\cdot, \cdot; \boldsymbol{\theta}_e)$ denotes
the density function of a bivariate copula of family $\mathcal{B}_e$
with parameters $\boldsymbol{\theta}_e$. The arguments $u_{i(e) \given
D(e)}$ and $u_{j(e) \given D(e)}$ of the bivariate copula density
functions $c_{\mathcal{B}_e}(\cdot, \cdot; \boldsymbol{\theta}_e)$ are
the values of the conditional copula cdfs,
\begin{align} \label{eq:u-data}
u_{i(e) \given D(e)} &:= C_{i(e) \given D(e)}\left(u_{i(e)};
T_{1{:}(k-1)}, \mathcal{B}_{1{:}(k-1)}, \boldsymbol{\theta}_{1{:}(k-1)}
\given\mathbf{u}_{D(e)}\right) \text{,} \nonumber\\
u_{j(e) \given D(e)} &:= C_{j(e) \given D(e)}\left(u_{j(e)};
T_{1{:}(k-1)}, \mathcal{B}_{1{:}(k-1)}, \boldsymbol{\theta}_{1{:}(k-1)}
\given\mathbf{u}_{D(e)}\right) \text{.}
\end{align}
The arguments $u_{i(e) \given D(e)}$ and $u_{j(e) \given D(e)}$ of a
pair copula $e \in E_k$ of level $k$ depend only on the specification
of the regular vine copula up to level $k-1$.
\end{definition}

\cite{cooke2001} showed that the regular vine copula density function
(see Definition \ref{def:rvine-cop}) is a valid $d$-variate probability
density function with uniform margins.

A regular vine copula is said to be {\it truncated at level $K$} if all
pair copulas conditional on $K$ or more variables are set to bivariate
independence copulas (\cite{brechmannczadoaas}). In that case, the pair
copula densities of trees $T_{K+1}, \ldots, T_{d-1}$ simplify to $1$
and do not affect the density of the truncated regular vine copula
anymore (see \eqref{eq:rvine-pdf}). This means that the first $K$
levels fully specify the truncated vine copula,
\[
\big(\mathcal{V}=(T_1,\ldots,T_K),\mathcal{B}_\mathcal{V} = (\mathcal
{B}_{1}, \ldots, \mathcal{B}_{K}), \boldsymbol{\theta}_{\mathcal{V}} =
(\boldsymbol{\theta}_{1}, \ldots, \boldsymbol{\theta}_{K})\big) \text{.}
\]

In general, regular vine copulas that differ in the tree structure or
in at least one pair copula family have different copula densities.
Notable exceptions from this rule include the multivariate Gaussian,
Student's t or Clayton copula, whose densities can be represented by
different pair copula constructions (cf. \cite{stoeber2012}).

\subsection{Likelihoods}
Definition \ref{def:rvine-cop} specifies regular vine copulas
tree-by-tree. Similarly, we introduce likelihoods of individual trees
and edges of a regular vine copula. The likelihoods are understood
given data $\mathbf{U} = (\mathbf{u}^1, \ldots, \mathbf{u}^T) \in
[0,1]^{T \times d}$. 
The likelihood of level $k$ depends only on the specification of the
regular vine copula up to level $k$. We write
\begin{align*}
L(\mathcal{B}_e, \boldsymbol{\theta}_e \given\mathbf{U}) &:= \prod
_{t=1{:}T} c_{\mathcal{B}_e}(u_{i(e) \given D(e)}^t, u_{j(e) \given
D(e)}^t; \boldsymbol{\theta}_e) \text{ and} \\
L(T_k, \mathcal{B}_{k}, \boldsymbol{\theta}_{k} \given\mathbf{U}) &:=
\prod_{e \in E_k} L(\mathcal{B}_e, \boldsymbol{\theta}_e \given\mathbf{U})
\end{align*}
for the likelihood of edge $e \in E_k$ and the likelihood of level $k$,
respectively. The pair $(u^t_{i(e) \given D(e)}, u^t_{j(e) \given
D(e)})$ denotes the $t$-th transformed observation (see \eqref
{eq:u-data}). The likelihood of a regular vine copula is then
calculated tree-by-tree
%
\begin{equation} \label{eq:rvine-lik}
L(\mathcal{V}, \mathcal{B}_{\mathcal{V}}, \boldsymbol{\theta}_{\mathcal
{V}} \given\mathbf{U}) = \prod_{T_k \in\mathcal{V}} L(T_k, \mathcal
{B}_{k}, \boldsymbol{\theta}_{k} \given\mathbf{U}) \text{.}
\end{equation}

\subsection{Pair Copula Families} \label{subsec:pair-copula-families}
\citet[Chapter 5]{joe2001} provides a collection of parametric copula
families for use in a pair copula construction. In our analyses, we
will consider the Independence (I), Gaussian (N) and Student's t (T)
copula as well as all rotations of the Gumbel (G) and Clayton (C)
copulas. Together, these pair copulas define the set of candidate
families $\mathbf{B}$. The transformations of the copulas' natural
parameters to their Kendall's $\tau$'s is provided in Tables 1 and 2 of
\cite{brechmannschepsmeier2013}.

Suppose a pair copula has a density function $c_0(u_1,u_2)$ and
strength of dependence parameter Kendall's $\tau= \tau_0$. The
90$^\circ$ rotation of a pair copula is defined by the density
$c_{90}(u_1,u_2) = c_0(1-u_1,u_2)$, the 180$^\circ$ rotation by
$c_{180}(u,v)=c_0(1-u_1,1-u_2)$ and the 270$^\circ$ rotation by
$c_{270}(u,v)=c_0(u_1,1-u_2)$. The 90$^\circ$ and 270$^\circ$ rotations
have Kendall's $\tau$'s $\tau_{90} = \tau_{270} = -\tau_0$, while the
Kendall's $\tau$ of the 180$^\circ$ rotation stays at $\tau_{180} = \tau_0$.


To shorten the notation in figures and tables, we will abbreviate the
pair copula families as above, possibly followed by the degrees of the
copulas' rotation and the values of their parameters in parentheses.
For example, C270($-0.8$) will indicate the 270$^\circ$ rotation of a
Clayton copula with strength of association parameter Kendall's $\tau= -0.8$.

\section[Bayesian Estimation of Regular Vine Copulas]{Bayesian
Estimation of Regular Vine Copulas Using Reversible Jump MCMC} \label
{ch:rvine-est}
Bayesian selection of regular vine copulas aims at estimating the joint
posterior distribution of the regular vine $\mathcal{V}$, pair copula
families $\mathcal{B}_{\mathcal{V}}$ and parameters $\boldsymbol{\theta
}_{\mathcal{V}}$.

The multi-layered composition of a regular vine copula and its density
function makes analytical inference infeasible. Instead, we use
reversible jump MCMC\vadjust{\goodbreak} (\cite{green1995}), which is an extension the
Metropolis--Hastings algorithm (\cite{metropolisetal1953,hastings1970})
to include the selection of models with different numbers of parameters
in the scope of inference, as a simulation-based approach to estimate
the posterior distribution. Convergence of the sampling chain to the
target distribution, here to the posterior distribution, is
theoretically established under regularity conditions.

\subsection{General Tree-by-Tree Model Selection} \label
{subsec:tree-by-tree-selection}
Our tree-by-tree model selection strategy first estimates the first
level of the regular vine copula, which consists of tree $T_1 = (N_1,
E_1)$ and the pair copula families $\mathcal{B}_{1}$ with parameters
$\boldsymbol{\theta}_{1}$. For each higher level $k = 2, \ldots, d-1$,
the density factorization $T_k = (N_k, E_k)$ and pair copula families
$\mathcal{B}_{k}$ with parameters $\boldsymbol{\theta}_{k}$ are
selected conditionally on the estimates of the lower levels $(T_1,
\mathcal{B}_1, \boldsymbol{\theta}_1)$ to $(T_{k-1}, \mathcal{B}_{k-1},
\boldsymbol{\theta}_{k-1})$, which remain unchanged from the previous steps.

\paragraph{Motivation of Tree-by-Tree Estimation} In the context of
model selection for regular vine copulas, sequential approaches exhibit
distinct strengths that make them more tractable than joint approaches.

Sequential approaches are much faster than joint approaches, as they
break the overall problem into a sequence of smaller problems that can
be solved more quickly. Table \ref{tb:no-of-rvines} shows the enormous
reduction of the regular vine search space, if a sequential procedure
is followed. The entries of Table \ref{tb:no-of-rvines} follow \cite
{morales2011}'s calculation of the number of vines and use the sum of
the number of spanning trees with $k$ nodes, $\sum_{k=2}^d k^{k-2}$, as
an upper bound of the size of the sequential search space. Here the
number of spanning trees is calculated using \cite{cayley1889}'s
formula. Furthermore, the reduced number of model alternatives improves
the convergence behavior of MCMC samplers as it allows for a quicker
exploration of the search space.
%
\begin{table}
\begin{center}
\begin{tabular}{c r r r r}
{} & \multicolumn{2}{c}{Vine Search Space} & \multicolumn{2}{c}{Vine
Copula Search Space} \\
{} & \multicolumn{1}{c}{Joint} & \multicolumn{1}{c}{Stepwise} &
\multicolumn{1}{c}{Joint} & \multicolumn{1}{c}{Stepwise} \\
\multicolumn{1}{l}{Dimension $d$} & \multicolumn{1}{c}{Selection} &
\multicolumn{1}{c}{Selection} & \multicolumn{1}{c}{Selection} &
\multicolumn{1}{c}{Selection} \\ \hline\hline
2 & 1 & 1 & 7 & 7 \\
3 & 3 & 3 & 1,029 & 154 \\
4 & 24 & $<$ 20 & 2,823,576 & $<$ 5,642 \\
5 & 480 & $<$ 145 & 1.3559e+11 & $<$ 305,767 \\
6 & 23,040 & $<$ 1,441 & 1.0938e+17 & $<$ 22,087,639 \\
7 & 2,580,480 & $<$ 18,248 & 1.4413e+24 & $<$ 1.9994e+9 \\
8 & 660,602,880 & $<$ 280,392 & 3.0387e+32 & $<$ 2.1789e+11 \\
9 & 3.8051e+11 & $<$ 5,063,361 & 1.0090e+42 & $<$ 2.7791e+13 \\
10 & 4.8705e+14 & $<$ 105,063,361 & 5.2118e+52 & $<$ 4.0632e+15
\end{tabular}
\caption{Size of the search space for vines $\mathcal{V}$ and vine
copulas $(\mathcal{V}, \mathcal{B}_{\mathcal{V}})$ with seven candidate
families, i.e., $|\mathbf{B}| = 7$, by dimension $d$.} \label{tb:no-of-rvines}
\end{center}
\end{table}

Furthermore, a tree-by-tree approach avoids a regular vine
copula-specific model identification issue. Different regular vine
copulas can be representatives of the same multivariate copula, the
most prominent example of which is the multivariate Gaussian copula
(\cite{cookekurowicka2006}). The tree-by-tree approach is characterized
by leaving previously selected trees unchanged and modifying only one
tree at a time. Under the tree-by-tree paradigm, there is only one
scenario in which the copula of the current state and proposed state
are the same with a non-zero probability: all pair copulas---those on
all previously selected trees and those on the current tree---are
either Gaussian or independent. These states can be easily detected and
collapsed into one state. 

\paragraph{Priors} Following our tree-by-tree estimation approach, the
priors are specified for each level $k=1,\ldots,d-1$. Given that the
proximity condition restricts which trees $T_k$ are allowed for a level
$k > 1$, these priors are inherently conditional on the selection on
the previous trees $T_1, \ldots, T_{k-1}$.

We choose a noninformative yet proper prior over the set $\mathbf
{STP}_k$ of all spanning trees that satisfy the proximity condition for
level $k$ for tree $T_k$, a sparsity-enforcing prior for the pair
copula families $\mathcal{B}_k$ and proper noninformative priors for
the parameters $\boldsymbol{\theta}_k$. We combine flat $(-1,1)$-priors
for the Kendall's $\tau$ parameters with flat $(0, \log(30))$-priors
for the logarithm of the degrees of freedom $\nu$ of Student's t pair copulas:
\begin{align*}
\pi(T_k) &\propto discrete \text{ } Uniform(\mathbf{STP}_k) \text{,} \\
\pi(\mathcal{B}_k \given T_k) &= \frac{\exp(-\lambda d_k)}{\sum
_{i=1}^{|E_k|} \sum_{d=0}^2 \exp(-\lambda d)} \propto\exp(-\lambda
d_k) \text{,} \\
\pi(\boldsymbol{\theta}_e \given T_k, \mathcal{B}_e) &\propto
\begin{cases}
Uniform_{(-1, 1)}(\tau_e) \text{ if $\mathcal{B}_e$ is a single
parameter copula,} \\
Uniform_{(-1, 1)}(\tau_e) \cdot\frac{\mathbf{1}_{(1,30)}(\nu_e) \cdot
\log(\nu_e)}{\int_1^{30} \log(x) dx} \text{ if $\mathcal{B}_e$ is the
Student's t copula,}
\end{cases}
\text{}
\end{align*}
where $d_k$ denotes the dimension of the parameter vector $\boldsymbol
{\theta}_k = (\boldsymbol{\theta}_{e;\mathcal{B}_e} \given e \in E_k)$
of the pair copula families $\mathcal{B}_k$ of level $k$. Analogously,
$d_e$ denotes the dimension of the parameter vector of the pair copula
family $\mathcal{B}_e$ of edge $e \in E_k$ and it holds that $d_k = \sum
_{e \in E_k} d_e$. Our prior on the pair copula families $\mathcal
{B}_k$ depends solely on the size $d_k$ of their parameter vectors
$\boldsymbol{\theta}_k$; if $\mathcal{B}_e$ is the independence copula,
it holds that $d_e = 0$.

The prior density $\pi$ of state $(T_k = (N_k, E_k), \mathcal{B}_k,
\boldsymbol{\theta}_k)$ results as
%
\begin{equation}
\pi(T_k, \mathcal{B}_k, \boldsymbol{\theta}_k) \propto\prod_{e \in
E_k} \exp(-\lambda d_e) \pi(\boldsymbol{\theta}_e \given T_k, \mathcal
{B}_e) \text{.}
\end{equation}

This prior gives the posterior distribution the following form:
\begin{align*} p(T_k, \mathcal{B}_k, \boldsymbol{\theta}_k \given
\mathbf{U}) &\propto\pi(\boldsymbol{\theta}_k \given T_k, \mathcal
{B}_k) \cdot\exp\left( \ell(T_k, \mathcal{B}_k, \boldsymbol{\theta}_k
\given\mathbf{U}) - \lambda d_k \right) \\ &\appropto\exp\left( \ell
(T_k, \mathcal{B}_k, \boldsymbol{\theta}_k \given\mathbf{U}) - \lambda
d_k \right) \text{,}
\end{align*}
where $\ell$ denotes the log-likelihood function and $\appropto$ means
``approximately proportional.'' At $\lambda= 0$, no shrinkage occurs
and the posterior mode estimate of level $k$ will approximate that
level's maximum likelihood estimate, while at $\lambda= 1$, the
posterior mode estimate of level $k$ will approximately minimize the
Akaike Information Criterion (AIC).

\paragraph{Posterior Distribution} The posterior distribution of level
$k$ given observed data $\mathbf{U}$ factorizes into the likelihood $L$
and prior density $\pi$:
\[
p(T_k, \mathcal{B}_k, \boldsymbol{\theta}_k \given\mathbf{U}) \propto
L(T_k, \mathcal{B}_k, \boldsymbol{\theta}_k \given\mathbf{U}) \cdot\pi
(T_k, \mathcal{B}_k, \boldsymbol{\theta}_k) \text{.}
\]
The tree-by-tree procedure requires the Bayesian posterior sample of
each tree to be collapsed into a single model estimate. We choose the
empirical mode of the sampled models $(T_k, \mathcal{B}_k)$ as the
model estimate, given that we chose our priors for their effects on the
posterior mode. The parameters are set to the means of the MCMC
posterior iterates of the selected model. Other centrality estimates
may be used as well.

\paragraph{Implementation} At each iteration $r=1,\ldots,R$, the
sampling mechanism performs a within-model move and a between-models
move. The within-model move updates all parameters $\boldsymbol{\theta
}_{1{:}k}$ of the regular vine copula, but leaves the pair copula
families $\mathcal{B}_{1{:}k}$ and tree structure $T_{1{:}k}$
unchanged. The between-models move operates only on level $k$ and
updates the tree structure $T_k$, pair copula families $\mathcal{B}_k$
along with the parameters $\boldsymbol{\theta}_k$.

The between-models move is implemented as a 50:50 mixture of two
mutually exclusive, collectively exhaustive (MECE) sub-routines: with a
50\% probability, a local between-models move updates only the pair
copula families $\mathcal{B}_k$ but leaves the tree structure $T_k$
unchanged (Algorithm \ref{alg:fam-update}). With the remaining 50\%
probability, a global between-models move updates the tree structure
$T_k$ along with the pair copula families~$\mathcal{B}_k$ (Algorithm
\ref{alg:tree-update}). Algorithm \ref{alg:fam-update} guarantees that
the proposal state differs in at least one pair copula family from the
current state; Algorithm \ref{alg:tree-update} guarantees that the
proposal state differs in at least one edge of tree $T_k$ from the
current state. This makes the proposals of the two sub-routines
mutually exclusive and gives the acceptance probability a tractable
analytical form that can be easily evaluated.

The between-models move is into two sub-routines because this allows an
intuitive interpretation of a local search (Algorithm \ref
{alg:fam-update}) and a global search (Algorithm \ref{alg:tree-update})
as well as optimizes the computational cost of these updates by
containing between-models moves that leave the tree structure unchanged
to a dedicated sub-routine.

\begin{algorithm}[Tree-by-Tree Bayesian Model Selection] \label
{alg:seq-outline} \quad
\begin{algorithmic}[1]
\FOR{{\bf each} level $k = 1, \ldots, d-1$}
\STATE Choose starting values: set tree $T_k = (N_k, E_k)$ to an
arbitrary tree that fulfills the proximity condition for level $k$; set
all pair copula families $\mathcal{B}_k$ of level $k$ to the
independence copula, i.e., $c_e(\cdot, \cdot) = 1$ for $e \in E_k$ and
set the parameter vector $\boldsymbol{\theta}_k$ of level $k$ to an
empty vector.
\FOR{{\bf each} MCMC iteration $r = 1,\ldots,R$}
\STATE Perform a within-model move: update all parameters $\boldsymbol
{\theta}_{1{:}k}$. Obtain $\boldsymbol{\theta}_{1{:}k}^{r,
\text{NEW}}$\break
through a Metropolis--Hastings step with random walk proposals:
\[
(T_k^r, \mathcal{B}_{k}^r, \boldsymbol{\theta}_{1{:}k}^r) = (T_k^{r-1},
\mathcal{B}_{k}^{r-1}, \boldsymbol{\theta}_{1{:}k}^{r, \text{NEW}})
\text{.}
\]

\STATE\label{alg:seq-outline:line-between-models-move} Perform a
between-models move: update the tree structure $T_k$ along with, or
only, the pair copula families $\mathcal{B}_{k}$ and parameters
$\boldsymbol{\theta}_{k}$ (Algorithms \ref{alg:fam-update}, \ref
{alg:tree-update}):
\[
(T_k^r, \mathcal{B}_{k}^r, \boldsymbol{\theta}_{k}^r) = (T_k^{r, \text
{NEW}}, \mathcal{B}_{k}^{r, \text{NEW}}, \boldsymbol{\theta}_{k}^{r,
\text{NEW}}) \text{.}
\]

\ENDFOR
\STATE Set the level $k$-estimate $(\hat{T}_k, \hat{\mathcal{B}}_{k},
\hat{\boldsymbol{\theta}}_{k})$ to the empirical mode of the posterior
sample $((T_k^r, \mathcal{B}_{k}^r, \boldsymbol{\theta}_{k}^r), r
= 1, \ldots, R)$:
\begin{itemize}
\item Set $\hat{T}_k$ and $\hat{\mathcal{B}}_{k}$ to the most
frequently sampled combination of $T_k$ and $\mathcal{B}_{k}$ in $
((T_k^r, \mathcal{B}_{k}^r), r = 1, \ldots, R)$.
\item Set $\hat{\boldsymbol{\theta}}_{k}$ to the sample mean of $
(\boldsymbol{\theta}_{k}^r, r \in\{1, \ldots, R\} \text{ with } T_k^r
= \hat{T}_k \text{ and}\break  \mathcal{B}_{k}^r = \hat{\mathcal{B}}_{k})$.
\end{itemize}
\STATE For all levels $l=1, \ldots, k-1$, update $\hat{\boldsymbol
{\theta}}_{l}$ and set it to the sample mean of \\ $(\boldsymbol
{\theta}_{l}^r, r \in\{1, \ldots, R\} \text{ with } T_k^r = \hat{T}_k
\text{ and } \mathcal{B}_{k}^r = \hat{\mathcal{B}}_{k})$.
\ENDFOR
\RETURN the stepwise Bayesian model estimate $(\hat{\mathcal{V}}, \hat
{\mathcal{B}}_{\mathcal{V}}, \hat{\boldsymbol{\theta}}_{\mathcal{V}})$,
where $\hat{\mathcal{V}} = (\hat{T}_1, \ldots, \hat{T}_{d-1})$, $\hat
{\mathcal{B}}_{\mathcal{V}} = (\hat{\mathcal{B}}_{1}, \ldots, \hat
{\mathcal{B}}_{d-1})$, and $\hat{\boldsymbol{\theta}}_{\mathcal{V}} =
(\hat{\boldsymbol{\theta}}_1, \ldots, \hat{\boldsymbol{\theta}}_{d-1})$.
\end{algorithmic}
\end{algorithm}

\subsection{Update of the Pair Copulas of Level $k$} \label{sec:fam-est}
This section describes a sub-routine of Algorithm \ref{alg:seq-outline}
to update the pair copula families $\mathcal{B}_{k}$ and parameters
$\boldsymbol{\theta}_{k}$ of level $k$ of a regular vine copula. This
updating step leaves the density factorization $\mathcal{V}$ unchanged.



This sub-routine first selects how many pair copulas will be updated
(Line \ref{alg:fam-update:N} of Algorithm~\ref{alg:fam-update}) and
then randomly selects which pair copulas will be updated---denoted by
$E \subseteq E_k$ in the remainder (Line \ref{alg:fam-update:E}). Next,
it generates a proposal that updates the selected pair copulas (Lines
\ref{alg:fam-update:prop-begin}--\ref{alg:fam-update:prop-end}), and,
lastly, accepts or rejects the proposal based on a Metropolis--Hastings
updating rules (Line \ref{alg:fam-update:MH-acceptance}).

The proposal step (Lines \ref{alg:fam-update:prop-begin}--\ref
{alg:fam-update:prop-end}) iterates through all selected pair copulas
$e \in E$. It first estimates the parameters $\boldsymbol{\theta}_{e;
\mathcal{B}_e^*}$ of each candidate pair copula family $\mathcal{B}_e^*
\in\mathbf{B} \setminus\mathcal{B}_e^r$, where the estimates are
denoted by $\hat{\boldsymbol{\theta}}_{e; \mathcal{B}_e^*}$. The
likelihoods of the different candidate copulas, $L(\mathcal{B}_e^*, \hat
{\boldsymbol{\theta}}_{e; \mathcal{B}_e^*} \given\mathbf{U})$, are
then used as the proposal probability weights of the respective copula
families: $q_\mathcal{B}(\mathcal{B}_e^r \to\mathcal{B}_e^*) \propto
L(\mathcal{B}_e^*, \hat{\boldsymbol{\theta}}_{e; \mathcal{B}_e^*}
\given\mathbf{U})$. After selecting a pair copula family, the proposal
parameters $\boldsymbol{\theta}_e^*$ are drawn from a normal
distribution centered at the parameter estimate $\hat{\boldsymbol{\theta
}}_{e; \mathcal{B}_e^*}$. The proposal distribution $q_N$ from which
$N$ is drawn (Line \ref{alg:fam-update:N}), the parameter estimation
procedure (Line \ref{alg:fam-update:fam-par-est}) and the covariance
matrix $\Sigma$ of the parameters' proposal distribution (Line \ref
{alg:fam-update:prop-par}) are MCMC tuning parameters.

Pair copula families that can model only positive or negative Kendall's
$\tau$'s such as the Clayton copula or Gumbel copula are extended to
cover the entire range $[-1, 1]$. This is implemented by replacing the
first argument $u_1$ of the copula density function $c(u_1, u_2)$ by
$1-u_1$ whenever the dependence parameter $\tau$ changes signs.

As this sub-routine and the one from Section \ref{sec:tree-est} produce
non-overlapping proposals, the acceptance probability follows as
%
\begin{equation} \label{eq:fam-AR}
\boldsymbol{\alpha} = \frac{L(T_k^r, \mathcal{B}_{k}^*, \boldsymbol
{\theta}_{k}^* \given\mathbf{U})}{L(T_k^r, \mathcal{B}_{k}^r,
\boldsymbol{\theta}_{k}^r \given\mathbf{U})} \cdot\frac{\pi(T_k^r,
\mathcal{B}_{k}^*, \boldsymbol{\theta}_{k}^*)}{\pi(T_k^r, \mathcal
{B}_{k}^r, \boldsymbol{\theta}_{k}^r)} \cdot\prod_{e \in E} \frac
{q_{\mathcal{B}}(\mathcal{B}_e^* \to\mathcal{B}_e^r) \cdot\phi_{(\hat
{\boldsymbol{\theta}}_{e;\mathcal{B}_e^r}, \Sigma)}(\boldsymbol{\theta
}_e^r)}{q_{\mathcal{B}}(\mathcal{B}_e^r \to\mathcal{B}_e^*) \cdot\phi
_{(\hat{\boldsymbol{\theta}}_{e;\mathcal{B}_e^*}, \Sigma)}(\boldsymbol
{\theta}_e^*)} \text{,}
\end{equation}
where $\phi_{\boldsymbol{\mu}, \Sigma}(\cdot)$ denotes the density
function of the truncated multivariate normal distribution with mean
$\boldsymbol{\mu}$ and covariance matrix $\Sigma$; the truncation is
assumed at the bounds of the respective parameters. Both the numerator
and denominator of the acceptance probability contain $q_N(N)$ as a
factor that cancels out and does not appear in \eqref{eq:fam-AR}, given
that the return move of any update must change the same number $N$ of
pair copulas as the outbound move.


\begin{algorithm}[Between-Models Move to Update the Pair Copula
Families $\mathcal{B}_{k}$ and Parameters $\boldsymbol{\theta}_{k}$]
\label{alg:fam-update} \quad
This is for the $r$th iteration of line \ref
{alg:seq-outline:line-between-models-move} of Algorithm \ref{alg:seq-outline}.
\begin{algorithmic}[1]
\STATE Select how many pair copulas are updated: $N \sim q_N(\cdot)$;
$N \in\{1, \ldots, |E_k|\}$. \label{alg:fam-update:N}
\STATE Select which pair copulas are updated: $E \subseteq E_k$ with
$|E| = N$. \label{alg:fam-update:E}
\FOR{{\bf each} selected pair copula $e \in E$} \label
{alg:fam-update:prop-begin}
\STATE\label{alg:fam-update:fam-par-est} For each candidate pair
copula family $\mathcal{B}_e \in\mathbf{B} \setminus\mathcal{B}_e^r$
estimate the copula parameter $\boldsymbol{\theta}_{e; \mathcal{B}_e}$
given the transformed data $(\mathbf{u}^{t=1{:}T}_{i(e) \given D(e)},
\mathbf{u}^{t=1{:}T}_{j(e) \given D(e)})$ and denote the parameter
estimate by $\hat{\boldsymbol{\theta}}_{e; \mathcal{B}_e}$.
\STATE Draw a new copula family $\mathcal{B}_e^* \in\mathbf{B}
\setminus\mathcal{B}_e^r$ from the proposal distribution
%
\begin{equation} \label{eq:fam-update:q-fam}
q_{\mathcal{B}}(\mathcal{B}_e^r \to\mathcal{B}_e^*) \propto L(\mathcal
{B}_e^*, \hat{\boldsymbol{\theta}}_{e; \mathcal{B}_e^*} \given\mathbf
{U}) \text{.}
\end{equation}

\STATE\label{alg:fam-update:prop-par} Draw new parameters $\boldsymbol
{\theta}_e^* \sim\mathcal{N}(\hat{\boldsymbol{\theta}}_{e; \mathcal
{B}_e^*}, \Sigma)$ from a normal distribution.
\STATE The proposal family for pair copula $e \in E$ is $\mathcal
{B}_e^*$ and the proposal parameter is $\boldsymbol{\theta}_e^*$.
\ENDFOR
\STATE\label{alg:fam-update:prop-end} The proposal families for level
$k$ are $\mathcal{B}_k^*$ and the proposal parameters are $\boldsymbol
{\theta}_k^*$, where
\begin{align*}
\mathcal{B}_k^* &= (\mathcal{B}_e^* \text{ for } e \in E \text{ and }
\mathcal{B}_e^r \text{ for } e \in E_k \setminus E), \\
\boldsymbol{\theta}_k^* &= (\boldsymbol{\theta}_e^* \text{ for } e \in
E \text{ and } \boldsymbol{\theta}_e^r \text{ for } e \in E_k \setminus
E) \text{.}
\end{align*}

\STATE Accept the proposal and set $(T_k^{r, \text{NEW}}, \mathcal
{B}_{k}^{r, \text{NEW}}, \boldsymbol{\theta}_{k}^{r, \text{NEW}}) =
(T_k^r, \mathcal{B}_{k}^*, \boldsymbol{\theta}_{k}^*)$ with probability
$\boldsymbol{\alpha}$ \xch{(}{(Eq.}\ref{eq:fam-AR}). If rejected, set $(T_k^{r,
\text{NEW}}, \mathcal{B}_{k}^{r, \text{NEW}}, \boldsymbol{\theta
}_{k}^{r, \text{NEW}}) = (T_k^r, \mathcal{B}_{k}^r, \boldsymbol{\theta
}_{k}^r)$. \label{alg:fam-update:MH-acceptance}
\end{algorithmic}
\end{algorithm}

\subsection{Joint Update of the Regular Vine and Pair Copulas of Level
$k$} \label{sec:tree-est}
This section presents a sub-routine of Algorithm \ref{alg:seq-outline}
to update the regular vine at level~$k$---that is, tree $T_k$---and the
pair copula families $\mathcal{B}_{k}$ and parameters $\boldsymbol
{\theta}_{k}$ of that level. Definition \ref{def:rvine} requires that
the lower level trees $T_1, \ldots, T_{k-1}$ of the regular vine be
specified before tree $T_k$ is estimated.



Algorithm \ref{alg:tree-update} describes our joint update procedure of
tree $T_k=(N_k, E_k)$ and the corresponding pair copula families
$\mathcal{B}_{k}$ and parameters $\boldsymbol{\theta}_{k}$. We denote
the set of all spanning trees with node set $N_k$ that satisfy the
proximity condition by $\mathbf{STP}_k$. The cardinality of this set is
computed using Kirchhoff's matrix tree theorem (\cite{kirchhoff1847})
to obtain the normalizing constants of the proposal and prior
distributions. In a first step, this sub-routine draws a new spanning
tree $T_k^* = (N_k, E_k^*) \in\mathbf{STP}_k \setminus T_k^r$ from the
proposal distribution $q_T(T_k^r \to T_k^*) \propto p^{|E_k^* \cap
E_k^r|} \cdot(1-p)^{|E_k^* \setminus E_k^r|}$ (Line \ref
{alg:tree-update:tree-prop}); this is just a random walk distribution
on the set of allowable regular vine trees of level $k$!\vadjust{\eject} Then, the
algorithm generates a proposal for the pair copula families $\mathcal
{B}_{k}^*$ and parameters $\boldsymbol{\theta}_{k}^*$ of this level as
in Algorithm \ref{alg:fam-update} (Lines \ref
{alg:fam-update:prop-begin}--\ref{alg:fam-update:prop-end} in Algorithm
\ref{alg:fam-update}; Lines \ref{alg:tree-update:fam-prop-begin}--\ref
{alg:tree-update:fam-prop-end} in Algorithm~\ref{alg:tree-update}). The
only difference is that all pair copula families in $\mathbf{B}$ are
permissible candidates here and the edges $e$ are different. We use the
notation $q_{\mathcal{B}}(\mathcal{B}_e^*)$ instead of $q_{\mathcal
{B}}(\mathcal{B}_e^r \to\mathcal{B}_e^*)$ to indicate the slightly
different proposal distributions. The entire proposal for level $k$ of
the regular vine copula consists of a new tree $T_k^*$, pair copula
families $\mathcal{B}_{k}^*$ and parameters $\boldsymbol{\theta
}_{k}^*$, and is accepted or rejected based on Metropolis--Hastings
updating rules (Line \ref{alg:tree-update:MH-acceptance}).

This sub-routine has three MCMC tuning parameters. The first is the
parameter $p$ of the proposal distribution for tree $T_k$: values $p >
0.5$ make tree proposals $T_k^*$ similar to the current tree $T_k^r$
more likely than proposals that are less similar to the current state.
The situation is reversed for values $p < 0.5$. The second tuning
parameter is the choice of the estimation procedure for the pair copula
parameter vectors (Line \ref{alg:tree-update:fam-par-est}) and the last
is the covariance matrix $\Sigma$ of the proposal distribution of the
parameters (Line~\ref{alg:tree-update:prop-par}).\looseness=1

The proposal mechanism of this update routine guarantees that the
proposed regular vine tree $T_k^*$ is different from the current state
$T_k$. This ensures that the proposals of this sub-routine and the one
of Section \ref{sec:fam-est} are mutually exclusive. Furthermore, the
proposal probability of the reverse move from tree $T_k^*$ to $T_k$ is
the same as the proposal probability of the away move, given that the
number of shared edges as well as differing edges is the same. As a
result, the acceptance probability of a proposal of this algorithm can
be easily obtained as
%
\begin{equation} \label{eq:tree-AR}
\mathbf{\alpha} = \frac{L(T_k^*, \mathcal{B}_{k}^*, \boldsymbol{\theta
}_{k}^* \given\mathbf{U})}{L(T_k^r, \mathcal{B}_{k}^r, \boldsymbol
{\theta}_{k}^r \given\mathbf{U})} \cdot\frac{\pi(T_k^*, \mathcal
{B}_{k}^*, \boldsymbol{\theta}_{k}^*)}{\pi(T_k^r, \mathcal{B}_{k}^r,
\boldsymbol{\theta}_{k}^r)} \cdot\frac{\prod_{e \in E_k^r} q_{\mathcal
{B}}(\mathcal{B}_e^r) \cdot\phi_{(\hat{\boldsymbol{\theta}}_{e;
\mathcal{B}_e^r}, \Sigma)}(\boldsymbol{\theta}_e^r)}{\prod_{e \in
E_k^*} q_{\mathcal{B}}(\mathcal{B}_e^*) \cdot\phi_{(\hat{\boldsymbol
{\theta}}_{e; \mathcal{B}_e^*}, \Sigma)}(\boldsymbol{\theta}_e^*)} \text{.}
\end{equation}

\begin{algorithm}[Between-Models Move for a Joint Update of Tree $T_k =
(N_k, E_k)$ and the Pair Copula Families $\mathcal{B}_{k}$ and
Parameters $\boldsymbol{\theta}_{k}$] \label{alg:tree-update} \quad\\
This is for the $r$th iteration of line \ref
{alg:seq-outline:line-between-models-move} of Algorithm \ref{alg:seq-outline}.
\begin{algorithmic}[1]
\STATE\label{alg:tree-update:tree-prop} Draw a new spanning tree
$T_k^* = (N_k, E_k^*) \in\mathbf{STP}_k \setminus T_k^r$ that
satisfies the proximity condition from the proposal distribution
%
\begin{equation} \label{eq:tree-update:q-tree}
q_T(T_k^r \to T_k^*) \propto p^{|E_k^* \cap E_k^r|} \cdot(1-p)^{|E_k^*
\setminus E_k^r|} \text{.}
\end{equation}

\FOR{each pair copula $e \in E_k^*$} \label{alg:tree-update:fam-prop-begin}
\STATE\label{alg:tree-update:fam-par-est} For each candidate pair
copula family $\mathcal{B}_e \in\mathbf{B}$ estimate the copula
parameter $\boldsymbol{\theta}_{e; \mathcal{B}_e}$ given the
transformed data $(\mathbf{u}^{t=1{:}T}_{i(e) \given D(e)}, \mathbf
{u}^{t=1{:}T}_{j(e) \given D(e)})$ and denote the parameter estimate by
$\hat{\boldsymbol{\theta}}_{e; \mathcal{B}_e}$.
\STATE Draw a new copula family $\mathcal{B}_e^* \in\mathbf{B}$ from
the proposal distribution
%
\begin{equation} \label{eq:tree-update:q-fam}
q_{\mathcal{B}}(\mathcal{B}_e^*) \propto L(\mathcal{B}_e^*, \hat
{\boldsymbol{\theta}}_{e; \mathcal{B}_e^*} \given\mathbf{U}) \text{.}
\end{equation}

\STATE\label{alg:tree-update:prop-par} Draw new parameters $\boldsymbol
{\theta}_e^* \sim\mathcal{N}(\hat{\boldsymbol{\theta}}_{e; \mathcal
{B}_e^*}, \Sigma)$ from a normal distribution.
\STATE The proposal family for pair copula $e \in E_k^*$ is $\mathcal
{B}_e^*$ and has proposal parameter $\boldsymbol{\theta}_e^*$.
\ENDFOR
\STATE\label{alg:tree-update:fam-prop-end} The proposal state is
$(T_k^*, \mathcal{B}_{k}^*, \boldsymbol{\theta}_{k}^*)$, where
\begin{equation*}
\mathcal{B}_k^* = (\mathcal{B}_e^* \given e \in E_k^*) \text{ and }
\boldsymbol{\theta}_k^* = (\boldsymbol{\theta}_e^* \given e \in E_k^*)
\text{.}
\end{equation*}

\STATE Accept the proposal and set $(T_k^{r, \text{NEW}}, \mathcal
{B}_{k}^{r, \text{NEW}}, \boldsymbol{\theta}_{k}^{r, \text{NEW}}) =
(T_k^*, \mathcal{B}_{k}^*, \boldsymbol{\theta}_{k}^*)$ with probability
$\boldsymbol{\alpha}$ \xch{(}{(Eq.}\ref{eq:tree-AR}). If rejected, set
$(T_k^{r, \text{NEW}}, \mathcal{B}_{k}^{r, \text{NEW}}, \boldsymbol
{\theta}_{k}^{r, \text{NEW}}) = (T_k^r, \mathcal{B}_{k}^r, \boldsymbol
{\theta}_{k}^r)$. \label{alg:tree-update:MH-acceptance}
\end{algorithmic}
\end{algorithm}

\subsection{Implementation in C++}
The model selection algorithms presented in this section are
implemented in a proprietary C++ software package. As the computational
cost of evaluating the likelihood and calculating parameter estimates
increases linearly with the number of observations in the data set,
these tasks are parallelized onto multiple CPU cores using OpenMP to
help reduce overall computing time. Furthermore, our software package
relies heavily on the tools and functionality provided by the boost
(\cite{boost}) and CppAD (\cite{cppad}) libraries: the boost library
contains a function that generates random spanning trees from a product
probability distribution based on edge weights, which we employ in our
implementation of Algorithm \ref{alg:tree-update}; the CppAD library
allows for automatic differentiation, which we use for parameter estimation.

\section{Simulation Study} \label{ch:simstudy}
We present a simulation study that compare our sequentially Bayesian
strategy with \cite{dissmann2010}'s frequentist model selection
algorithm, the independence model and the maximum likelihood estimate
(MLE) of the multivariate Gaussian copula. The comparisons with the
independence model and Gaussian copula illustrate that vine copulas are
relevant dependence models that significantly improve model fit over
simpler standard models, while the comparison with Di{\ss}mann's vine
copula estimates highlights the improved model selection capabilities
of our method.

Our simulation study uses copula data generated from four different
six-dimensional vine copulas (Table \ref{tb:simstudy-true-models} of
Appendix \ref{ch:app-sim}). These scenarios cover a wide range of
dependence structures: the first two cover general cases of
multivariate dependence, while the third and fourth scenario are
special cases to investigate detailed characteristics of our model
selection method. Scenario 3 consists of only one level, which means
that all variables are conditionally independent. It also means that
the true model lies in the search space of the first level of our
selection procedure, so that this scenario can be used to validate our
proposed scheme empirically. Scenario 4 is has only Gaussian pair
copulas, which makes it a vine copula-representation of the
multivariate Gaussian copula. As a result, this scenario allows for an
isolated evaluation of the pair copula family selection aspect of our
method, given that the multivariate Gaussian copula results from any
vine density factorization $\mathcal{V}$ as long as all pair copula
families are Gaussian.


We generate 100 data sets consisting of 500 independent samples from
the respective copula distribution of each scenario and allow the pair
copula families listed in Section~\ref{subsec:pair-copula-families} as
candidates.

\subsection{Choice of the Benchmark Algorithm}
The models selected by \cite{dissmann2010}'s algorithm serve as a
benchmark. This algorithm follows a stepwise frequentist approach that
selects each tree $T_k$, $k = 1,\ldots,5$ as the maximum spanning tree
using absolute values of Kendall's $\tau$ of the variable pairs as edge
weights. The pair copula families are selected to optimize the AIC
copula-by-copula and the parameters are set to their individual maximum
likelihood estimates.

Di{\ss}mann's algorithm and ours share their tree-by-tree selection
strategy. However, there are two major differences between our
approaches: firstly, Di{\ss}mann follows a heuristic scheme to select
the tree structure $\mathcal{V}$, while we follow a proper Bayesian
selection scheme on each level $k$; secondly, Di{\ss}mann selects the
pair copula families on an edge-by-edge basis, whereas we place priors
on the distribution of the pair copula families across an entire level
$k$ to simultaneously select of all edges of that level.

\subsection{Configuration of Our Reversible Jump MCMC Sampler} \label
{ss:simstudy-configuration}
We use the shrinkage prior introduced in Section \ref
{subsec:tree-by-tree-selection} with shrinkage parameter $\lambda= 1$.
The posterior mode estimates of each level $k$ will then be
approximately AIC-optimal.

We use our reversible jump MCMC Algorithm \ref{alg:seq-outline} from
Section \ref{ch:rvine-est} to generate $R = 50{,}000$ posterior samples
for each level $k = 1,\ldots,5$ of the 6-dimensional regular vine
copula. The MCMC tuning parameters are summarized in Table \ref
{tb:tuning-parameters} of Appendix \ref{ch:app-sim}. Furthermore, we
apply a re-weighting on the proposal probabilities (see \eqref
{eq:fam-update:q-fam} and \eqref{eq:tree-update:q-fam}) of the pair
copula families in the sub-routines of Algorithms \ref{alg:fam-update}
and \ref{alg:tree-update} to improve the mixing behavior of the
sampling chain. This is achieved by ensuring that the ratio of smallest
and biggest the proposal probabilities is bounded from below by 0.05,
\[
\frac{\min_{\mathcal{B}^* \in\mathbf{B} \setminus\mathcal{B}^r}
q_{\mathcal{B}}(\mathcal{B}^r \to\mathcal{B}^*)}{\max_{\mathcal{B}^*
\in\mathbf{B} \setminus\mathcal{B}^r} q_{\mathcal{B}}(\mathcal{B}^r
\to\mathcal{B}^*)} \geq0.05\ \text{ and }\  \frac{\min_{\mathcal{B}^*
\in\mathbf{B}} q_{\mathcal{B}}(\mathcal{B}^*)}{\max_{\mathcal{B}^* \in
\mathbf{B}} q_{\mathcal{B}}(\mathcal{B}^*)} \geq0.05 \text{,
respectively.}
\]

\subsection{Evaluation of the Results}
\begin{figure}[p!]
\includegraphics{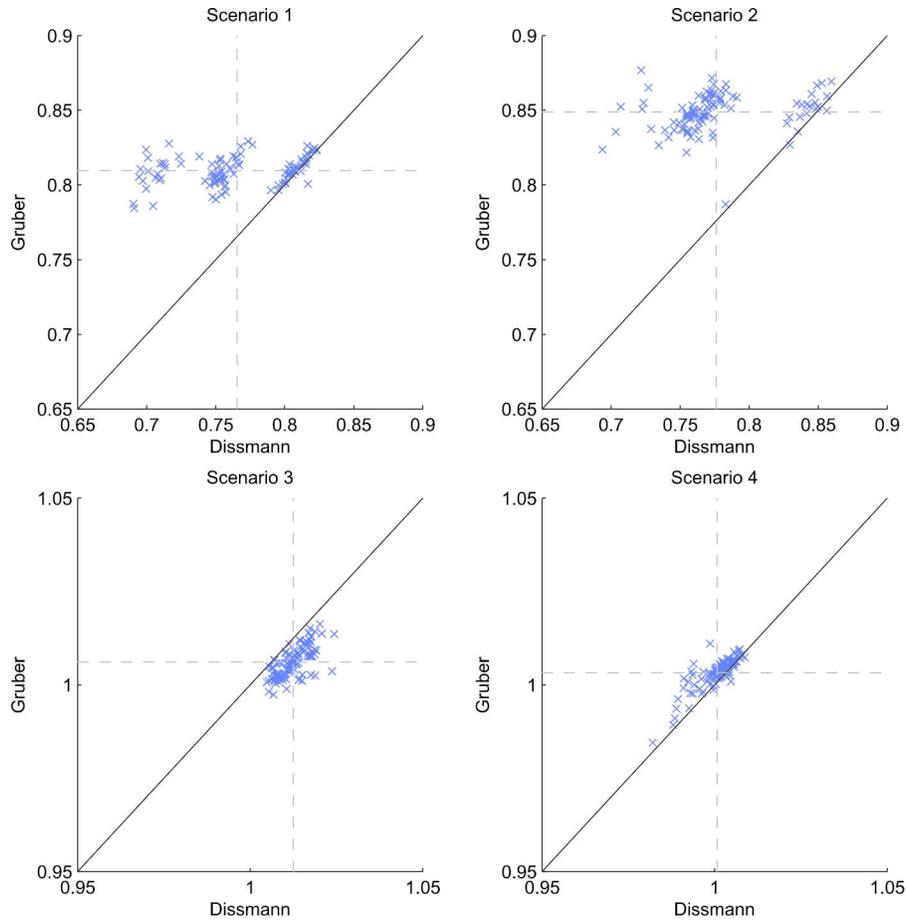}
\caption{Comparison of relative log-likelihoods of our method and
Di{\ss}mann's. The dashed lines indicate the respective averages across
all 100 replications.} \label{fig:simstudy-loglik-comparison}
\end{figure}

The results are based on 100 replications of the estimation procedures
with independently generated data sets of size 500 from the four
scenarios and are summarized in Table \ref
{tb:simstudy-results-summary}. The fitting capabilities of our
algorithm and \cite{dissmann2010}'s are measured by the log-likelihood
of the estimated models. Knowing the underlying ``true'' models, we can
also calculate the ratio of the estimated log-likelihoods and the true
log-likelihoods to evaluate how well the selection methods perform in
absolute terms. Figure \ref{fig:simstudy-loglik-comparison} compares
the performance of our Bayesian strategy with Di{\ss}mann's heuristic:
markers above the diagonal line indicate replications in which our
Bayesian model estimate has a higher likelihood than Di{\ss}mann's.
%
%
\begin{table}[p!]
\tabcolsep=5.8pt
\begin{center}
\begin{tabular}{l c c c c}
{} & Scenario 1 & Scenario 2 & Scenario 3 & Scenario 4 \\ \hline\hline
Gruber $>$ Di{\ss}mann (out of 100)& 97 & 98 & 0 & 86 \\ \hline
Gruber rel. loglik (in \%) & {\bf81.0} & {\bf84.9} &100.6 & {\bf
100.3}\\
Di{\ss}mann rel. loglik (in \%) & 76.6 & 77.6 & {\bf101.3}& 100.1 \\
Gaussian MLE rel. loglik (in \%) & 64.0 & 67.6 & 84.7 & 100.1 \\
Independence rel. loglik (in \%) & 0 & 0 & 0 & 0 \\[6pt]
\end{tabular}
\caption{Number of replications in which our algorithm's estimate has a
higher likelihood than Di{\ss}mann's; average percentage of the true
log-likelihood of the estimated vine copulas, the multivariate Gaussian
copula and the independence model.} \label{tb:simstudy-results-summary}
\end{center}
\end{table}

\paragraph{Scenarios 1 and 2} The log-likelihoods of the models
selected by our algorithm average 81\% and 85\% of the log-likelihoods
of the true models, Di{\ss}mann's model estimates average 77\% and 78\%
, and the multivariate Gaussian copula averages 64\% and 68\%,
respectively (Table \ref{tb:simstudy-results-summary}). While neither
Di{\ss}mann's method nor ours selects the correct tree $T_1$ in any
replication, these numbers still show that model fit is improved
significantly by our approach, and that Di{\ss}mann's estimates are
already more suitable models than the multivariate Gaussian copula.

Figure \ref{fig:simstudy-loglik-comparison} shows that within each
scenario, the relative log-likelihoods of our model selection procedure
are distributed much more narrowly about the mean compared with Di{\ss
}mann's. Furthermore, in 97 (Scenario 1) and 98 (Scenario 2) out of 100
replications, our method's estimates perform better than Di{\ss}mann's
(Table \ref{tb:simstudy-results-summary}). Together, this shows that
our model selection strategy is more robust and performs consistently
better.

\paragraph{Scenario 3} The regular vine copula of Scenario 3 is
truncated to the first level. As a result, the true model lies in the
search space of the first step of our tree-by-tree model selection
procedure, which makes this a test case to validate our implementation.

Our shrinkage prior effectively avoids over-fitting, given that, on
average, only 0.8 of the 10 pair copulas on levels $k=2{:}5$ are
selected as non-independence copulas, and in 43 out of 100
replications, all pair copulas on levels $k=2{:}5$ are selected as
independence copulas. Furthermore, in all 100 replications, the
posterior mode estimate has the true model's tree structure $T_1$.
Di{\ss}mann's procedure is more prone to over-fitting with, on average,
3.3 out of the 10 pair copulas on levels $k=2{:}5$ being
non-independence copulas and only 2 out of 100 replications selecting
all pair copulas as independence copulas.

The log-likelihoods of the estimated models by our algorithm average
101\% as do the log-likelihoods from Di{\ss}mann's models. This is an
excellent result that confirms the validity of our model selection
scheme and shows that it is implemented correctly. The consistently
slightly higher log-likelihoods of Di{\ss}mann's model estimates are
based on over-fitting. This scenario confirms our method as superior to
Di{\ss}mann's, as it is important for an effective selection method to
identify sparse patterns. The multivariate Gaussian copula lags behind
with an average relative log-likelihood of 85\% even though it is the
model that has the most parameters.

{\em Detailed Analysis of the MCMC Output of Level 1 of Replication 1.}
After discarding the first 2,500 iterations as burn-in, the posterior
mode model has a posterior probability of 58\% (Model 29; Table \ref
{tb:v1-l1-rep1-hist}) and all posterior samples, after burn-in, have
the correct tree structure $T_1$. The selected pair copula families
agree with the correct pair copula families, except for the family of
edge $3,6$: Model 29 selects a Gaussian pair copula, Model 30 selects
the 180 degree rotation of the Gumbel copula, and Model 31 selects the
Student's t copula. The fact that the correct model, Model 31, has only
4\% posterior probability can be attributed to our shrinkage prior,
given that the log-likelihoods of these three models are nearly
identical. Figure \ref{fig:v1-l1-posterior-sample} of Appendix \ref
{ch:app-sim} illustrates the MCMC mixing behavior using the model index
and log-likelihood trace plots.

\begin{table}
\begin{center}
\begin{tabular}{l r r r}
Posterior Model & \multicolumn{1}{c}{\bf29} & \multicolumn{1}{c}{30} &
\multicolumn{1}{c}{31} \\ \hline\hline
Posterior probability (in \%) & {\bf58.4} & 32.1 & 4.0 \\
Average relative log-likelihood (in \%) & {\bf100.2} & 100.1 & 100.4 \\
Number of parameters & {\bf5} & 5 & 6 \\
Correct tree $T_1$ & {\bf yes} & yes & yes \\
\end{tabular}
\caption{Scenario 3, Replication 1. Histogram table of selected models
with an empirical posterior probability of at least 1\%. The posterior
probabilities and relative log-likelihoods are quoted in percentage
points.} \label{tb:v1-l1-rep1-hist}
\end{center}
\end{table}

\paragraph{Scenario 4} Both model selection procedures select models
that average about 100\% of the log-likelihoods of the true model. This
extraordinary performance can be explained by a peculiarity of vine
copulas: if all pair copulas are Gaussian or independence copulas, the
vine copula equals a multivariate Gaussian copula irrespective of the
density factorization $\mathcal{V}$. As a result, the selection of the
density factorization $\mathcal{V}$ does not play a role in selecting
suitable vine copula models here. Our sequential Bayesian procedure
selects, on average, 13.8 out of the 15 pair copulas as either Gaussian
or independence copulas, while Di{\ss}mann's procedures comes in second
at 12.7 out of 15. This result, together with the high relative
log-likelihoods, suggests that both algorithms perform similarly well
at selecting suitable pair copula families.

\paragraph{Conclusion} Both algorithms perform equally well in fitting
a vine copula to Gaussian data. Our tree-level Bayesian approach
improves model selection of general regular vine copulas, which are not
independent of the selected tree structure $\mathcal{V}$. The large
performance gap between Scenarios 1, 2 and the special cases of
Scenarios 3, 4 shows the limits of our tree-by-tree approach towards
the selection of the tree structure $\mathcal{V}$. We acknowledge that
our model selection scheme does not yet represent the definitive answer
to the model selection challenge. Nevertheless, our proposed selection
scheme consistently selects better-fitting models than existing
selection strategies and is better at detecting sparsity patterns for
model reduction than \cite{dissmann2010}'s frequentist method.

\subsection{Analysis of the Computational Complexity and Runtime}
\paragraph{Computational Complexity} The computation of a single model
estimate by our algorithm using the set-up described in Section \ref
{ss:simstudy-configuration} consists of 50,000 MCMC updates for each
level of the vine copula. These sum up to 250,000 MCMC updates for the
five levels of a six-dimensional regular vine copula. Each MCMC update
consists of a between-models move and a within-model move: the
between-models move consists of estimating the parameters and
calculating the likelihood of each pair copula and each candidate pair
copula family, and an additional evaluation of the likelihood after
drawing a proposal parameter; the within-model move brings another
evaluation of the likelihood of each pair copula.

\paragraph{Computing Facilities and Runtime} The simulation study was
performed on a Linux cluster with AMD Opteron-based 32-way nodes using
2.6 GHz dual core CPUs for parallel processing. The Linux cluster is
hosted by the Leibniz-Rechenzentrum der Bayerischen Akademie der
Wissenschaften near Munich, Germany. It took approximately 10 hours to
execute our stepwise Bayesian selection strategy for a six-dimensional
data set of size 500, while it took 5--6 seconds to execute \cite
{dissmann2010}'s heuristic and less than 0.1 seconds to estimate the
correlation matrix of the multivariate Gaussian copula. It may be noted
that the runtime of our procedure could be cut significantly by
reducing the number of MCMC iterations. Our analyses suggest that
convergence is achieved quickly so $R\in[15{,}000, 30{,}000]$ will be
adequate choices in practice.

\paragraph{Recommendations for Researchers} In most studies,
researchers have to strike a balance between getting quick, or getting
more accurate results. With that in mind, we propose the following
approach to decide which dependence model to use. In a quick first
analysis, estimate a multivariate Gaussian copula and select a regular
vine copula using \cite{dissmann2010}'s heuristic, which can be
completed within a few seconds. If the log-likelihoods of both models
are similar, use the multivariate Gaussian copula as a ``good enough''
standard model (see Scenario 4). However, if the log-likelihood of the
selected vine copula is substantially higher than the one of the
Gaussian copula, perform a sequential Bayesian analysis using our
method for more accurate and more robust results, and better sparsity
detection (see Scenarios 1--3).

\section{Example: Portfolio Asset Returns} \label{ch:case}
We consider a diversified portfolio that invests in multiple asset
classes using iShares exchange-traded funds (ETFs) and commodity
trusts. The daily log-returns of each investment are modeled by a
univariate time series. The joint multivariate characteristics are
modeled by a regular vine copula and a multivariate Gaussian copula.

We learn the copulas using one year's worth of data and then use the
selected copulas together with the marginal time series to obtain joint
multivariate step-ahead forecasts for six months. The quality of the
forecasts is measured by comparing the forecast accuracy of various
portfolio metrics with the actual realizations.

\subsection{Description of the Data}
The data set contains adjusted daily closing prices of nine iShares
ETFs, $j=1{:}9$, and covers the time period from January 2013 through
June 2014.\footnote{The data were downloaded from \url
{http://finance.yahoo.com}.} The training set consists of 252
observations from January through December 2013 ($t=1{:}252$); the test
set consists of 124 observations from January through June 2014
($t=253{:}376$). The nine ETFs form a well-diversified portfolio that
invests in multiple asset classes and can be easily replicated by
retail investors (Table \ref{tb:etf-overview}). Three of the funds
invest in U.S. equities ($j=1,2,3$), two funds in U.S. treasuries
($j=4,5$), two funds in U.S. real estate through real estate investment
trusts (REITs, $j=6,7$), and two funds are commodity trusts investing
in gold and silver ($j=8,9$).\footnote{More details on the selected
funds can be found on the iShares homepage at \texttt{\surl
{http://www.\\ishares.com/us/index}}.}

\begin{table}
\begin{center}
\begin{tabular}{c c p{2.05in} p{1.8in}}
$j$ & Symbol & Name & Exposure \\ \hline\hline
$1$ & IVV & iShares Core S\&P 500 ETF & Large-cap U.S. stocks \\
$2$ & IJH & iShares Core S\&P Mid-Cap ETF & Mid-cap U.S. stocks \\
$3$ & IJR & iShares Core S\&P Small-Cap ETF & Small-cap U.S. stocks \\
\hline
$4$ & HYG & iShares iBoxx \$ High Yield Corporate Bond ETF & High yield
corporate bonds \\
$5$ & LQD & iShares iBoxx \$ Investment Grade Corporate Bond ETF & U.S.
investment grade corporate bonds \\ \hline
$6$ & RTL & iShares Retail Real Estate Capped ETF & U.S. retail
property real estate stocks and REITs \\
$7$ & REZ & iShares Residential REIT Capped ETF & U.S. residential real
estate stocks and REITs \\ \hline
$8$ & SLV & iShares Silver Trust & Silver \\
$9$ & IAU & iShares Gold Trust & Gold
\end{tabular}
\caption{Overview of the ETFs selected for the real data example. The
exposure information is taken from the iShares homepage.} \label
{tb:etf-overview}
\end{center}
\end{table}
%

\subsection{Marginal Time Series}
We model the daily log-returns $y_{j,t}$, $t=1,2,\ldots$ of each series
$j=1{:}9$ using a variance discounting dynamic linear model (DLM; \cite
[Chapter 10.8]{West1997}). The DLM is a fully Bayesian time series
model that has closed-form posterior and forecast distributions, and
the parameters are learned on-line. The following updating equations
are adapted from Table 10.4 of \cite{West1997}.

\paragraph{Model Formulation} The general DLM models each time series
$y_{j,t}$, $j=1{:}9$, by
\begin{align}
y_{j,t} &= \mathbf{F}_{j,t}' \boldsymbol{\theta}_{j,t} + \nu_{j,t}, &
\nu_{j,t} &\sim N(0, \lambda_{j,t}^{-1}), \label{eq:dlm-obs} \\
\boldsymbol{\theta}_{j,t} &= \mathbf{G}_{j,t} \boldsymbol{\theta
}_{j,t-1} + \boldsymbol{\omega}_{j,t}, & \boldsymbol{\omega}_{j,t}
&\sim N(\mathbf{0}, \mathbf{W}_{j,t}), \label{eq:dlm-system} \\
\lambda_{j,t} &= \lambda_{j,t-1} \frac{\eta_{j,t}}{\beta_j}, & \eta
_{j,t} &\sim Beta\left(\frac{\beta_j n_{j,t-1}}{2}, \frac{(1-\beta_j)
n_{j,t-1}}{2}\right), \label{eq:dlm-precision}
\end{align}
with observation equation (\ref{eq:dlm-obs}). Equations \eqref
{eq:dlm-system} and \eqref{eq:dlm-precision} describe the evolutions of
the states $\boldsymbol{\theta}_{j,t}$, a vector with $p_j$ entries,
and $\lambda_{j,t}$, a positive scalar, where the innovations $\nu
_{j,t}$, $\boldsymbol{\omega}_{j,t}$ and $\eta_{j,t}$ are assumed
mutually independent and independent over time. The predictors $\mathbf
{F}_{j,t}$ are a vectors of size $p_j$ and the state evolution matrices
$\mathbf{G}_{j,t}$ are of dimensions $p_j \times p_j$. The parameters
$\mathbf{W}_{j,t}$, $n_{j,t-1}$ and $\beta_j$ of the state evolutions
\eqref{eq:dlm-system}--\eqref{eq:dlm-precision} are explained in the
next paragraph.

\paragraph{Forward Filtering} The information set at time $t$ is
denoted by $\mathcal{D}_t$. Suppose that at time $t-1$, a normal--gamma
prior for $(\boldsymbol{\theta}_{j,t}, \lambda_{j,t})$, given
information $\mathcal{D}_{t-1}$ has density
%
\begin{equation} \label{eq:dlm-prior}
\pi_{j,t}(\boldsymbol{\theta}_{j,t}, \lambda_{j,t}) := N(\boldsymbol
{\theta}_{j,t} \given\mathbf{a}_{j,t}, \mathbf{R}_{j,t}/(c_{j,t}
\lambda_{j,t})) \cdot G(\lambda_{j,t} \given r_{j,t}/2, r_{j,t} c_{j,t}
/ 2) \text{,}
\end{equation}
and parameters $\mathbf{a}_{j,t} \in\mathbb{R}^{p_j}$, $\mathbf
{R}_{j,t} \in\mathbb{R}^{p_j \times p_j}$, $r_{j,t} > 0$ and $c_{j,t}
> 0$; at $t=0$, the initial prior parameters are $\mathbf{a}_{j,1}$,
$\mathbf{R}_{j,1}$, $r_{j,1}$ and $c_{j,1}$. As $y_{j,t}$ is observed
at time $t$, the information set is updated to $\mathcal{D}_t$ and the
posterior distribution of $(\boldsymbol{\theta}_{j,t}, \lambda_{j,t})$,
given information $\mathcal{D}_{t}$ follows as a normal--gamma
%
\begin{equation} \label{eq:dlm-posterior}
p_{j,t}(\boldsymbol{\theta}_{j,t}, \lambda_{j,t}) := N(\boldsymbol
{\theta}_{j,t} \given\mathbf{m}_{j,t}, \mathbf{C}_{j,t}/(s_{j,t}
\lambda_{j,t})) G(\lambda_{j,t} \given n_{j,t}/2, n_{j,t} s_{j,t} / 2)
\end{equation}
with parameters $\mathbf{m}_{j,t} = \mathbf{a}_{j,t} + \mathbf{A}_{j,t}
e_{j,t} \in\mathbb{R}^{p_j}$, $\mathbf{C}_{j,t} = (\mathbf{R}_{j,t} -
\mathbf{A}_{j,t} \mathbf{A}_{j,t}' q_{j,t}) z_{j,t} \in\mathbb{R}^{p_j
\times p_j}$, $n_{j,t} = r_{j,t} + 1 > 0$ and $s_{j,t} = z_{j,t}
c_{j,t} > 0$, where $e_{j,t} = y_{j,t} - \mathbf{F}_{j,t}'\mathbf
{a}_{j,t} \in\mathbb{R}$ is the forecast error, $q_{j,t} = c_{j,t} +
\mathbf{F}_{j,t}' \mathbf{R}_{j,t} \mathbf{F}_{j,t} > 0$ is the
forecast variance factor, $\mathbf{A}_{j,t} = \mathbf{R}_{j,t} \mathbf
{F}_{j,t} / q_{j,t} \in\mathbb{R}^{p_j}$ is the adaptive coefficient
vector, and $z_{j,t} = (r_{j,t} + e_{j,t}^2 / q_{j,t}) / n_{j,t} > 0$
is the volatility update factor (see Table 10.4 of \cite{West1997}).
The step-ahead priors $(\boldsymbol{\theta}_{j,t+1}, \lambda_{j,t+1}
\given\mathcal{D}_t)$ at time $t$ follow from the system equations
\eqref{eq:dlm-system}--\eqref{eq:dlm-precision} as evolutions of the
posterior states $(\boldsymbol{\theta}_{j,t}, \lambda_{j,t} \given
\mathcal{D}_t)$. The normal--gamma step-ahead prior density $\pi
_{j,t+1}$ is as in \eqref{eq:dlm-prior} with $t$ evolved to $t+1$ and parameters
$r_{j,t+1} = \beta_j n_{j,t}$, $c_{j,t+1} = s_{j,t}$, $\mathbf
{a}_{j,t+1} = \mathbf{G}_{j,t+1} \mathbf{m}_{j,t}$, $\mathbf{R}_{j,t+1}
= \mathbf{G}_{j,t+1} \mathbf{C}_{j,t} \mathbf{G}_{j,t+1}' + \mathbf
{W}_{j,t+1}$ and $\mathbf{W}_{j,t+1} = \frac{1-\delta_j}{\delta_j}
\mathbf{G}_{j,t+1} \mathbf{C}_{j,t} \mathbf{G}_{j,t+1}'$. The discount
factors $\beta_j, \delta_j \in(0,1)$ inflate the prior variances in
the state evolution steps and determine the model's responsiveness to
new observations.

\paragraph{Forecasting} The forecast distribution of $y_{j,t+1}$ at
time $t$ and with information $\mathcal{D}_t$ follows as a
non-standardized Student's t distribution $T_{\text{non std}}(\nu, \mu,
\sigma^2)$ with $\nu= r_{j,t+1}$ degrees of freedom, mean $\mu=
\mathbf{F}_{j,t+1}'\mathbf{a}_{t+1}$ and variance $\sigma^2=\mathbf
{F}_{j,t+1}' \mathbf{R}_{j,t+1} \mathbf{F}_{j,t+1} + c_{j,t+1}$ by
integration of the observation equation (\ref{eq:dlm-obs}) over the
prior distributions of the states $(\boldsymbol{\theta}_{j,t+1}, \lambda
_{j,t+1})$. The non-standardized t distribution is a location--scale
transformation
%
\begin{equation} \label{eq:dlm-forecast}
T_{\text{non std}}(\nu, \mu, \sigma^2) = \mu+ \sqrt{\sigma^2} \cdot
T_{\nu}
\end{equation}
of a t distribution $T_{\nu}$ with $\nu$ degrees of freedom. In the
remainder, we will denote the forecast distribution of $y_{j,t+1}$ at
time $t$ by $T_{j,t+1}$.

\begin{figure}
\includegraphics{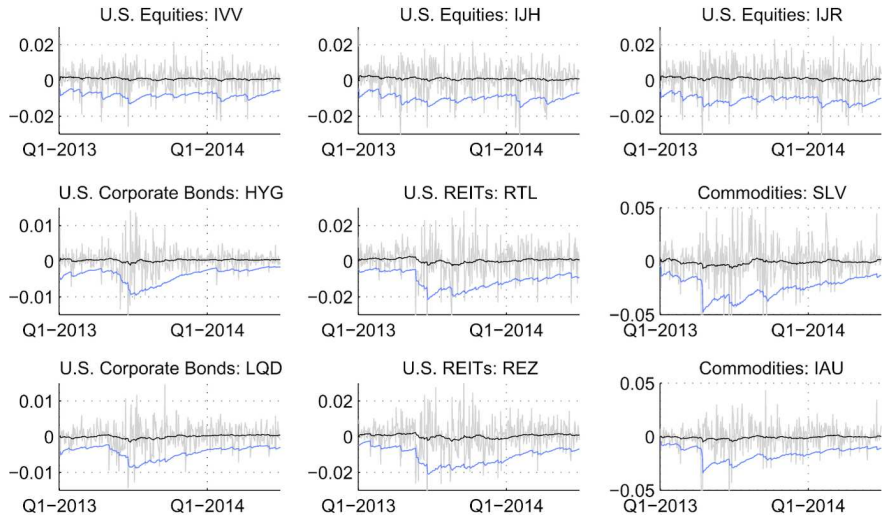}
\caption{Realized returns (gray), forecast means (black) and 10\% value
at risk (blue).} \label{fig:DLMs}
\end{figure}

\paragraph{Model Choice} We use a local-level DLM that assumes $\mathbf
{F}_{j,t} = 1$ and has random walk evolutions $\mathbf{G}_{j,t}=1$ for
all $j$ and $t$. The discount factors are set to $\beta_j = 0.96$ and
$\delta_j = 0.975$ for all $j$ to balance responsiveness to new
observations with sufficient robustness for reliable forecasts. We
start the analysis with the initial prior parameters of each series
$j=1{:}9$ set to $\mathbf{a}_{j,1} = 0$, $\mathbf{R}_{j,1}=10^{-6}$,
$r_{j,1}=10$, $c_{j,1}=10^{-5}$. Figure~\ref{fig:DLMs} shows the
sequential step-ahead forecasts and realized daily log-returns as well
as the 10\% quantiles of the forecast distributions as the daily value
at risk of each series $j=1{:}9$. 

\subsection{Estimation of the Dependence Models}
Copula modeling is a two-step process: first, marginal models remove
within-series effects from the data $y_{j,t}$ to obtain i.i.d.---within
each series $j$---uniform noise $u_{j,t} := T_{j,t}(y_{j,t})$; second,
a copula is selected to describe across-series dependence effects of
the multivariate transformed $U(0,1)$ data $\mathbf{u}_t = (u_{1,t},
\ldots, u_{9,t})'$, $t=1,2,\ldots$. 

\paragraph{Sequential Bayesian Selection} We use our model selection
scheme to estimate a 9-dimensional regular vine copula using the
$t=1{:}252$ observations from 2013.\vadjust{\goodbreak} We apply the same priors and
configuration of the sample as in Section \ref{ch:simstudy}. The
selected model is shown in Table \ref{tb:case-vine-est} of Appendix \ref
{ch:app-case}: it has 10 Gaussian pair copulas, 0 Student's t copulas,
8 Gumbel copulas, 0 Clayton copulas, and 18 Independence copulas. The
pair copulas of levels $k \geq5$ are selected as Independence copulas
and omitted in Table \ref{tb:case-vine-est}.

\paragraph{Di{\ss}mann's Frequentist Selection} We compare against our
Bayesian tree-by-tree strategy to \cite{dissmann2010}'s frequentist
heuristic as we did in the simulation study of Section \ref
{ch:simstudy}. Di{\ss}mann's vine copula is noticeably less
parsimonious with only 13 Independence pair copulas and 5 Gaussian
copulas, 7 Student's t copulas, 7 Gumbel copulas, and 4 Clayton
copulas. The selected model is shown in Table \ref
{tb:case-dissmann-est} of Appendix \ref{ch:app-case}.

\paragraph{Multivariate Gaussian Copula} For reference, we also
included a maximum likelihood estimate of the multivariate Gaussian
copula in our comparison. The estimated correlation matrix is shown in
Table \ref{tb:case-gauss-est} of Appendix \ref{ch:app-case}.

\subsection{Analysis of Portfolio Forecasts}

\paragraph{Sampling from the Joint Forecast Distribution} Samples
$\mathbf{\hat{y}}_t^{n=1{:}N} = (\hat{y}_{1,t}^n, \ldots,\break  \hat
{y}_{9,t}^n)'$ from the joint forecast distribution are generated by
transforming samples $\mathbf{u}^n=(u_1^n, \ldots, u_9^n)'$ from the
copula to the observation scale through inverse cdfs of the marginal
forecast distributions, $\hat{y}_{j,t}^n := T_{j,t}^{-1}(u_j^n)$.

\paragraph{In-Sample Analysis} Consider a portfolio that invests
equally in the ETFs from Table~\ref{tb:etf-overview} and the weights
$\mathbf{w}_t := (w_{1,t}, \ldots, w_{9,t})' = (\frac{1}{9}, \ldots,
\frac{1}{9})'$ are maintained throughout the time period $t=1{:}252$.
For each $t$, we draw $N=10{,}000$ samples $\mathbf{\hat
{y}}_t^{n=1{:}N}$ from the joint forecast distribution to simulate the
portfolio returns $\hat{r}_t^n := \mathbf{w}_t' \mathbf{\hat{y}}_t^n$.
We compute the 10\% value at risk as the 10\% sample quantile and the
expected portfolio return as the sample mean of $\hat{r}_t^{n=1{:}N}$.
This allows an evaluation of the adequacy of the joint multivariate
model, which consists of the nine marginal DLMs and the selected copula.

\begin{figure}
\includegraphics{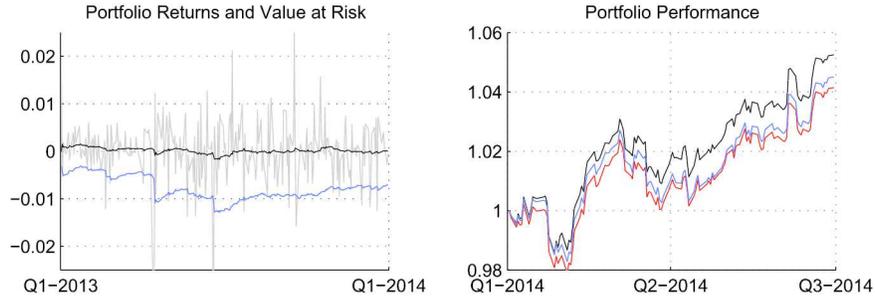}
\caption{(Left) Realized portfolio returns (gray) vs. forecast
means (black) and 10\% value at risk (blue) from our sequential
Bayesian vine copula model during the training date range $t=1{:}252$.
(Right) Portfolio performance of investment strategy
(\ref{eq:max-sharpe-ratio}) under our Bayesian vine copula model (black;
highest line), Di{\ss}mann's vine model (red; lowest line) and the
Gaussian model (blue; middle line) during the test date range
$t=253{:}376$.} \label{fig:portfolio-forecasts}
\end{figure}

Figure \ref{fig:portfolio-forecasts} (left) shows that the predicted
quantities from our sequential Bayesian vine copula model are in line
with the actual portfolio returns. The actual portfolio return is under
the predicted 10\% quantile of the simulated portfolio return
distribution that uses our sequential Bayesian vine copula on 8.7\%, or
22 out of 252 days; if Di{\ss}mann's vine copula or the multivariate
Gaussian copula are used in conjunction with the same marginal DLMs,
the actual portfolio return is under the predicted 10\% quantile on
8.3\%, or 21 out of 252 days; if the independence copula is used, the
10\% value at risk is exceeded 20\%, or 50 out of 252 days.
%

\paragraph{Out-of-Sample Analysis} During the out-of-sample period
from January through June 2014, $t=253:376$, we investigate the
performance of a dynamic portfolio whose weights $\mathbf{w}_t$ are
updated daily to maximize the predicted Sharpe ratio (\cite{sharpe1966}):
%
\begin{equation} \label{eq:max-sharpe-ratio}
\max_{\mathbf{w}_t} \widehat{SR}_t(\mathbf{w}_t) \text{ subject to }
\sum_{j=1{:}9} w_{j,t} = 1 \text{ and } w_{j,t} \in(0.05, 0.25) \text{
for all } j=1{:}9 \text{.}
\end{equation}
Here $\widehat{SR}_t(\mathbf{w}_t)$ (see \eqref{eq:def-sharpe-ratio})
is the estimate of the annualized Sharpe ratio of a portfolio with
investment weights $\mathbf{w}_t$, where $\boldsymbol{\hat{\mu}}_t$ is
the sample mean and $\widehat{\Sigma}_t$ is the sample covariance
matrix of the simulated joint forecasts $\mathbf{\hat{y}}_t^{n=1{:}N}$.
Again, $N=10{,}000$ samples were used:
%
\begin{equation} \label{eq:def-sharpe-ratio}
\widehat{SR}_t(\mathbf{w}_t) :
= \frac{252 \cdot\mathbf{w}_t' \boldsymbol{\hat{\mu}}_t}{\sqrt{252
\cdot\mathbf{w}_t' \widehat{\Sigma}_t \mathbf{w}_t}} \text{.}
\end{equation}
Our calculation of the Sharpe ratio in \eqref{eq:def-sharpe-ratio} is
under the assumption of a zero-return risk-free asset. The optimization
constraints $w_{j,t} \in(0.05, 0.25)$ in \eqref{eq:max-sharpe-ratio}
refer to minimum and maximum weights of individual assets and are
typical restrictions that aim at protecting investors from undue
accumulation of risk.

When the regular vine copula selected by our sequential Bayesian
procedure is used as the joint model's dependence model, the realized
annualized Sharpe ratio of the investment strategy (\ref
{eq:max-sharpe-ratio}) during the out-of-sample period $t=253{:}376$ is
${SR}=1.95$; if the multivariate Gaussian copula is used to inform the
investment decisions, the realized Sharpe ratio is ${SR}=1.67$; if
Di{\ss}mann's frequentist vine model is used, the realized Sharpe ratio
is ${SR}=1.53$. In addition, the realized nominal return of the
portfolio driven by our sequential Bayesian vine copula is higher than
the returns of the portfolios using the Gaussian copula or Di{\ss
}mann's copula (Figure \ref{fig:portfolio-forecasts} (right)). This
example provides additional evidence of our sequential Bayesian vine
model as the most reliable model for use in a real-life context. 

\subsection{Considerations on Use in Practice}
%
The computing time of updating the univariate DLMs is a few
milliseconds, which is negligible in the context of daily portfolio
risk rebalancing. In contrast, the computational burden of estimating
the regular vine copula is much higher---in our example, our Bayesian
model selection could take as long as up to a day to complete.

We suggest that the dependence model be update in monthly, quarterly,
or semi-annually intervals only. Under the assumption that the
dependence structure of financial asset returns is only slowly
changing, this is a prudent way to proceed. Even though the long
computation time of our Bayesian strategy might be a deterrent to
implementing our approach, the benefits of increasing a portfolio's
performance in nominal as well as risk-adjusted terms will quickly pay
for the investment in computing time.

The combination of univariate DLMs with a regular vine copula as the
dependence model provides a robust framework for forecasting, yet is
highly responsive to new observations. This can be seen, for example,
in the way the value at risk changes instantly on the day of a large
negative market move. Furthermore, it is a distinct strength of regular
vine copulas to be able to model asymmetric dependence characteristics
along with various tail dependence characteristics in one model. Our
example shows that a regular vine copula-driven model can help in
decision making to achieve superior investment performance as well as
improved risk forecasts.

\section{Concluding Remarks} \label{ch:conclusions}
We discussed a Bayesian approach to model selection of regular vine
copulas and presented a reversible jump Markov chain Monte Carlo-based
algorithm to facilitate posterior sampling. A key feature of our
approach, sequential model selection in the levels $k$ reduces the
search space for candidate models to a fraction of the search space for
joint selection to keep the computational run time at an acceptable level.

A simulation study (Section \ref{ch:simstudy}) demonstrated that our
Bayesian model selection approach is superior to \cite{dissmann2010}'s
frequentist one. The better performance of our Bayesian selection
scheme can be attributed to its simultaneous and prior-informed
selection of the pair copula families $\mathcal{B}_k$ of a given level
$k$, while Di{\ss}mann's algorithm selects them one-by-one. 
In addition to the simulation study, a real data example (Section \ref
{ch:case}) illustrated how regular vine copulas can be used to achieve
superior portfolio risk forecasts and investment decisions.

Our estimation procedure extends previously available inference methods
for regular vine copulas in two significant ways. Our Bayesian
tree-by-tree strategy allows the selection of the pair copula families
$\mathcal{B}_{\mathcal{V}}$ from an arbitrary set of candidate families
$\mathbf{B}$, which is a non-trivial extension of \cite
{smithminalmeidaczado2009}'s indicator-based approach that can only
detect (conditional) pairwise independencies. Furthermore, we present
the first Bayesian inference method for selecting the regular vine
$\mathcal{V}$ as the building plan of the pair copula construction
jointly with the pair copula families $\mathcal{B}_{\mathcal{V}}$. A
major selling point of our approach is that we demonstrated its
superiority to existing procedures in a simulation study under
controlled conditions (see Section \ref{ch:simstudy}) as well as in an
application study using real data (see Section \ref{ch:case}).

Sequential model selection schemes can fail to select the correct
model. This is illustrated, e.g., in Section \ref{ch:simstudy} by the
failure of the selected models to have relative log-likelihoods close
to 100\% in Scenarios 1 and 2. Current research aims at developing a
fully Bayesian model selection scheme to estimate all levels of a
regular vine copula jointly as well as allowing for time-varying
dependence effects.

\appendix


\FloatBarrier
\section{Supplements to the Simulation Study} \label{ch:app-sim}
%
\begin{table}[!h]
\begin{center}
\scalebox{.92}{
\begin{tabular}{l p{4.25in}}
Algorithm & Tuning Parameters \\ \hline\hline
\ref{alg:fam-update} & $q_N(N = k) = \frac{1}{3.5} \log\left( 1 - \frac
{1-e^{-3.5}}{|E_k| e^{-3.5} + k (1 - e^{-3.5})} \right)$, where $|E_k|$
denotes the number of pair copulas of level $k$\\
\ref{alg:fam-update}, \ref{alg:tree-update} & Parameter estimation is
done by matching the Kendall's $\tau$ parameter to the sample Kendall's
$\tau$. The degrees of freedom $\nu$ of a Student's t pair copula is
maximum likelihood estimated on a discrete grid. \\
\ref{alg:fam-update}, \ref{alg:tree-update} & $\Sigma= 0.0125^2$ for
the Kendall's $\tau$ of single-parameter copulas;\newline$\Sigma=
\begin{pmatrix} 0.0125^2 & 0 \\ 0 & 0.1^2
\end{pmatrix}
$ for the $(\tau, \log\nu)$ parameter vector of the Student's t copula
\\
\ref{alg:tree-update} & $p =0.667$ \\
\end{tabular}
}
\caption{MCMC tuning parameters used in the simulation study and real
data example.} \label{tb:tuning-parameters}
\end{center}
\end{table}

\begin{table}[t!]
\begin{center}
\scalebox{.92}{
\begin{tabular}{l l l l}
\multicolumn{1}{l}{Scenario 1} & \multicolumn{1}{l}{Scenario 2} &
\multicolumn{1}{l}{Scenario 3} & \multicolumn{1}{l}{Scenario 4} \\
\hline\hline
$c_{1,2}$ N(0.59) & $c_{1,2}$ T(0.54, 5) & $c_{1,2}$ N(0.41) &
$c_{1,2}$ N(0.41) \\
$c_{2,3}$ C(0.71) & $c_{1,3}$ C90(--0.67) & $c_{2,3}$ C(0.50) &
$c_{2,3}$ N(0.49) \\
$c_{3,4}$ C180(0.80) & $c_{1,4}$ C180(0.64) & $c_{3,4}$ C180(0.50) &
$c_{2,4}$ N(--0.33) \\
$c_{3,5}$ N(--0.71) & $c_{1,5}$ N(--0.59) & $c_{3,5}$ N(--0.33) &
$c_{3,5}$ N(--0.26) \\
$c_{3,6}$ T(0.65, 3) & $c_{1,6}$ T(0.54, 6) & $c_{3,6}$ T(0.49, 5) &
$c_{3,6}$ N(0.13) \\
$c_{1,3 \given2}$ G(0.75) & $c_{2,3 \given1}$ G(0.71) & & $c_{1,3
\given2}$ N(0.59) \\
$c_{2,4 \given3}$ N(0.41) & $c_{2,4 \given1}$ G270(--0.71) & &
$c_{2,5 \given3}$ N(0.13) \\
$c_{2,5 \given3}$ C270(--0.60) & $c_{2,5 \given1}$ C270(--0.60) & &
$c_{3,4 \given2}$ N(0.41) \\
$c_{2,6 \given3}$ N(--0.37) & $c_{2,6 \given1}$ N(--0.45) & & $c_{5,6
\given3}$ N(--0.33) \\
$c_{1,4 \given2,3}$ T(0.26, 5) & $c_{3,4 \given1,2}$ T(0.30, 8) & &
$c_{1,5 \given2,3}$ N(0.26) \\
$c_{1,5 \given2,3}$ N(--0.26) & $c_{3,5 \given1,2}$ N(--0.30) & &
$c_{2,6 \given3,5}$ N(--0.41) \\
$c_{1,6 \given2,3}$ C90(--0.56) & $c_{3,6 \given1,2}$ C90(--0.43) & &
$c_{4,5 \given2,3}$ N(0.19) \\
$c_{4,6 \given1,2,3}$ N(0.13) & $c_{4,5 \given1,2,3}$ N(0.19) & &
$c_{1,6 \given2,3,5}$ N(0.49) \\
$c_{5,6 \given1,2,3}$ C(0.20) & $c_{4,6 \given1,2,3}$ C(0.43) & &
$c_{4,6 \given2,3,5}$ N(0.41) \\
$c_{4,5 \given1,2,3,6}$ G180(0.52) & $c_{5,6 \given1,2,3,4}$
G180(0.50) & & $c_{1,4 \given2,3,5,6}$ N(--0.33) \\ \hline
17 parameters & 18 parameters & 6 parameters & 15 parameters
\end{tabular}
}
\caption{The vine copulas used in the simulation study. The parameters
shown are the Kendall's $\tau$ and the degrees of freedom $\nu$.} \label
{tb:simstudy-true-models}
\end{center}
\end{table}

\begin{figure}[!t]
\includegraphics{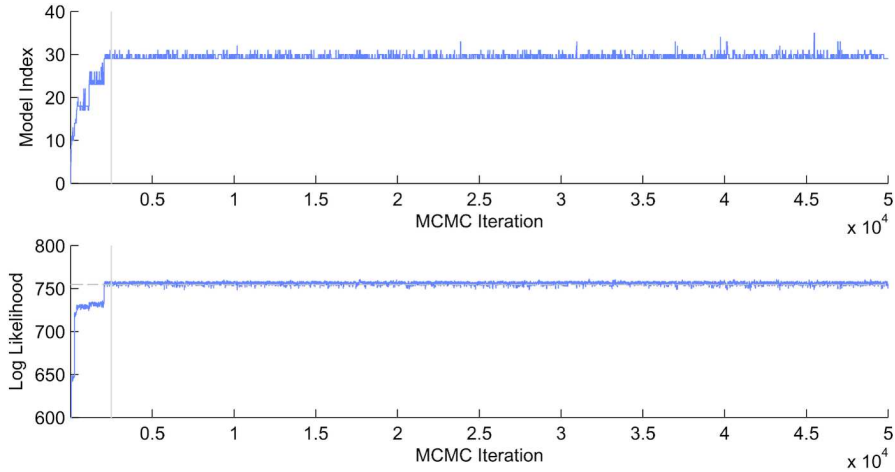}
\caption{Model index and log-likelihood trace plot of Replication 1 of
Scenario 3. The horizontal line in the lower plot indicates the true
model's log-likelihood; the vertical lines show the burn-in period.}
\label{fig:v1-l1-posterior-sample}\vspace*{50pt}
\end{figure}

\clearpage
\section{Supplements to the Data Example} \label{ch:app-case}

\begin{table}[!h]
\begin{center}
\scalebox{.92}{
\begin{tabular}{l l l l}
Tree $T_1$ & Tree $T_2$ & Tree $T_3$ & Tree $T_4$\\ \hline\hline
$c_{1,2}$ N(0.74) & $c_{1,3;2}$ N(0.07) & $c_{1,9;5,7}$ I &
$c_{1,8;6,7,9}$ G180(0.06)\\
$c_{1,4}$ N(0.46) & $c_{1,5;7}$ G270(--0.10) & $c_{2,6;1,7}$ G180(0.08)
& $c_{2,9;1,6,7}$ I \\
$c_{1,7}$ N(0.39) & $c_{1,6;7}$ N(0.14) & $c_{3,7;1,2}$ I &
$c_{3,6;1,2,7}$ I \\
$c_{2,3}$ N(0.79) & $c_{2,7;1}$ G(0.12) & $c_{4,5;1,7}$ N(0.27) &
$c_{4,6;1,5,7}$ I \\
$c_{5,7}$ G(0.19) & $c_{4,7;1}$ G180(0.13) & $c_{5,6;1,7}$ N(0.12) &
$c_{5,9;1,6,7}$ G(0.12) \\
$c_{6,7}$ N(0.51) & $c_{6,8;9}$ I & $c_{7,8;6,9}$ I \\
$c_{6,8}$ G(0.09) & $c_{7,9;6}$ I \\
$c_{8,9}$ N(0.71) \\
\end{tabular}
}
\caption{Sequential Bayesian estimate of the regular vine copula, given
the training data $t=1{:}252$. This tables shows the Kendall's $\tau$
parameters of the pair copulas.} \label{tb:case-vine-est}
\end{center}
\end{table}

\begin{table}[!h]
\begin{center}
\scalebox{.92}{
\begin{tabular}{l l l l}
Tree $T_1$ & Tree $T_2$ & Tree $T_3$ \\ \hline\hline
$c_{1,2}$ N(0.74) & $c_{1,5;4}$ T(--0.18, 12.8) & $c_{1,3;2,7}$ N(0.06)
\\
$c_{1,4}$ N(0.46) & $c_{1,7;2}$ G180(0.08) & $c_{1,9;4,5}$ G180(0.03) \\
$c_{2,3}$ T(0.79, 8.67) & $c_{2,4;1}$ I & $c_{2,5;1,4}$ I \\
$c_{2,7}$ T(0.41, 8.58) & $c_{2,6;7}$ N(0.16) & $c_{3,6;2,7}$ I \\
$c_{4,5}$ T(0.39, 4.93) & $c_{3,7;2}$ G270(--0.05) & $c_{4,7;1,2}$
G180(0.13) \\
$c_{5,9}$ G(0.16) & $c_{4,9;5}$ C(0.05) & $c_{4,8;5,9}$ C180(0.05) \\
$c_{6,7}$ T(0.51, 14.5) & $c_{5,8;9}$ C270(--0.05) \\
$c_{8,9}$ N(0.70) \\
\\
Tree $T_4$ & Tree $T_5$ & Tree $T_6$ \\ \hline\hline
$c_{1,6;2,3,7}$ I & $c_{2,8;1,4,5,9}$ T(--0.04, 16.7) &
$c_{3,9;1,2,4,5,7}$ I \\
$c_{1,8;4,5,9}$ G(0.06) & $c_{2,5;1,2,4,7}$ I & $c_{5,6;1,2,3,4,7}$
T(0.10, 14.3) \\
$c_{2,9;1,4,5}$ C(0.04) & $c_{4,6;1,2,3,7}$ I & $c_{7,8;1,2,4,5,9}$ I \\
$c_{3,4;1,2,7}$ I & $c_{7,9;1,2,4,5}$ I \\
$c_{5,7;1,2,4}$ G(0.17)
\end{tabular}
}
\caption{Regular vine copula selected by Di{\ss}mann's heuristic, given
the training data $t=1{:}252$. This tables shows the Kendall's $\tau$
parameters of the pair copulas.} \label{tb:case-dissmann-est}
\end{center}
\end{table}

\begin{table}[!h]
\begin{center}
\scalebox{.92}{
\begin{tabular}{c}
$\Sigma=
\begin{pmatrix}
\mathbf{1} & 0.92 & 0.89 & 0.64 & 0.02 & 0.58 & 0.61 & 0.21 & 0.16 \\
0.92 & \mathbf{1} & 0.95 & 0.61 & 0.03 & 0.59 & 0.61 & 0.21 & 0.17 \\
0.89 & 0.95 & \mathbf{1} & 0.58 & 0.00 & 0.54 & 0.57 & 0.20 & 0.15 \\
0.64 & 0.61 & 0.58 & \mathbf{1} & 0.37 & 0.53 & 0.57 & 0.23 & 0.23 \\
0.02 & 0.03 & 0.00 & 0.37 & \mathbf{1} & 0.32 & 0.32 & 0.14 & 0.20 \\
0.58 & 0.59 & 0.54 & 0.53 & 0.32 & \mathbf{1} & 0.74 & 0.20 & 0.18 \\
0.61 & 0.61 & 0.57 & 0.57 & 0.32 & 0.74 & \mathbf{1} & 0.18 & 0.16 \\
0.21 & 0.21 & 0.20 & 0.23 & 0.14 & 0.20 & 0.18 & \mathbf{1} & 0.90 \\
0.16 & 0.17 & 0.15 & 0.23 & 0.20 & 0.18 & 0.16 & 0.90 & \mathbf{1}
\end{pmatrix}
$
\end{tabular}
}
\caption{Maximum likelihood estimate of the correlation matrix of the
multivariate Gaussian copula, given the training data $t=1{:}252$.}
\label{tb:case-gauss-est}
\end{center}
\end{table}


\begin{acknowledgement}
The numerical computations were performed on a Linux cluster supported
by DFG grant INST 95/919-1 FUGG.

We thank the Editor, Associate Editor and three unnamed reviewers for
constructive comments on our original manuscript.
\end{acknowledgement}

\end{document}